\documentclass[11pt]{article}
\usepackage{fullpage}
\usepackage{amsmath,amsthm,amsfonts,amssymb,dsfont,bm}
\usepackage{enumerate,color,xcolor}
\usepackage[colorlinks,linkcolor=cyan,citecolor=blue]{hyperref}
\usepackage{url}

\numberwithin{equation}{section}
\theoremstyle{plain}

\newtheorem{theorem}{Theorem}

\newtheorem{lemma}[theorem]{Lemma}
\newtheorem{proposition}[theorem]{Proposition}

\newtheorem{remark}[theorem]{Remark}

\begin{document}
\begin{center}  
{\bf\Large Sequencing, task failures, and capacity 
when failures are driven by a non-homogeneous Poisson process}
\end{center}

\author{}
\begin{center}
  Lingjiong Zhu\,\footnote{Department of Mathematics, Florida State University, Tallahassee, Florida, United States of America; zhu@math.fsu.edu
 }
   Anand Paul\,\footnote{Department of Information Systems and Operations Management, University of Florida, Gainesville, Florida, United States of America; paulaa@ufl.edu
 }
   Haldun Aytug\,\footnote{Department of Mathematics, Department of Information Systems and Operations Management, University of Florida, Gainesville, Florida, United States of America; aytugh@ufl.edu
 }
\end{center}

\begin{center}
 \today
\end{center}

\begin{abstract}
We study the optimal sequencing of a batch of tasks on a machine subject to random disruptions driven by a non-homogeneous Poisson process (NHPP), such that every disruption requires the interrupted task to be re-processed from scratch, and partially completed work on a disrupted task is wasted. The NHPP models random disruptions whose frequency varies systematically with time. In general the time taken to process a given batch of tasks depends on the order in which the tasks are processed.  We find conditions under which  the simplest possible sequencing rules - shortest processing time first (SPT) and longest processing time first (LPT) - suffice to minimize the completion time of a batch of tasks.  
\end{abstract}



\section{Problem Description and Literature Review}
\label{sec-Related work}
We study the problem of sequencing a fixed set of tasks on a single machine with the objective of minimizing the expected completion time of the entire set of tasks, when the machine is subject to breakdowns at points in time that are generated by a nonhomogeneous Poisson process (NHPP). The machine processes tasks one by one in a fixed sequence, proceeding from one task to the next one in the sequence after completely processing the preceding task. A task finishes processing only when the machine processes it from start to finish without a single disruption. If the machine breaks down in the course of processing a task, the time spent in partially processing it is accounted for in the total processing time, and the processing of the task starts over from scratch; the same protocol is followed if the machine should break down more than once in the course of processing any job. In general, each iterated processing of the same task is governed by a different probability distribution of completion time since the NHPP is a nonstationary process; in contrast, if breakdowns occurred at points in time generated by a Poisson process,  the iterated processing of a task is an exact probabilistic replica of the initial process. The objective is to determine the sequence in which the tasks should be processed in order to minimize the expected total completion time of the entire set of tasks. Because we are only interested in the effects of sequencing on batch completion times (makespan), we do not model down times for repair. Past work on this problem may be categorized  based on whether disruptions
occur according to a Poisson process, a renewal process, or a NHPP. We note that
technically, for a renewal process to be stationary, the time to the first
disruption/renewal is required to have the equilibrium distribution associated with the renewal
distribution.

A key aspect of the problem is that if the machine is interrupted in the course of processing a task, then  the task will have to be reprocessed from scratch and the time invested in partially processing the task is wasted.   This model of processing is called {\it preempt-repeat} in the scheduling literature (\cite{Birge1990}). An entirely different model ({\it preempt-resume}) arises if we assume that the processor can resume processing an interrupted task at the point where it left off (see \cite{righter94} for a survey of research in the preempt-resume model). 
\begin{remark}
To connect the preempt-repeat model with practice, we remark that in the consulting experience of one of the authors, an overhaul process for aircraft required an expensive CNC machine that was used for the production of a titanium alloy aerospace part; the machine suffered breakdowns down in the course of part production and whenever it did, a fresh part production cycle had to be restarted after the machine was restored to operating condition since it was not feasible  to resume part production starting from a partially produced part. For another application, we note that  data sent over a communications network routinely fail to reach its destinations. This behavior is expected and data network protocols, such as TCP/IP, are written to accommodate failures. Failures in transmission occur due to buffer overflows in router hardware and noise in communication lines. The typical scenario is that a data packet has to be resent in its entirety when it gets corrupted or lost during transmission, and this furnishes another example of the preempt-repeat model.  
\end{remark}


The preempt-repeat model is notoriously hard to analyze when the times between disruptions are independent and identically distributed (i.i.d.) but not exponential. Work on optimal sequencing in this case
includes \cite{Adiri1989}, \cite{Kasap2006, Kasap2008} and \cite{aytugandpaul2012}. Adiri et al. \cite{Adiri1989} consider the case of a single disruption, and
Kasap et al. \cite{Kasap2006, Kasap2008} proved, for two tasks, that when the times between disruptions has an
increasing density function, LPT minimizes the expected makespan. \cite{aytugandpaul2012} show that SPT minimizes expected makespan and flow-time of an arbitrary number of tasks when the uptime distribution is a mixture of exponentials, a subclass of the decreasing failure rate (DFR) distributions. Sequencing tasks with random processing times and disruptions according to a
Poisson process have been studied extensively \cite{Frostig1991, Cai2003, Cai2004, Birge1990}. 
A different, but interesting, set of questions in the preempt-repeat regime is explored in \cite{asmussenetal2008}. They consider the case when the uninterrupted task processing times are random and derive tight asymptotic relations between the total completion time and the sum of uninterrupted processing times. They prove that if the task time distribution is unbounded, then the distribution of total completion time is heavy-tailed. 

One typically asks if there are operational policies that can minimize wasted time when a processor faced with multiple tasks fails, with failures generated by a random---but stationary---process. \cite{Adiri1989} show that if the failures occur at points in time generated by a renewal process such that the time between failures is distributed according to a convex (concave) distribution function and a batch of tasks has known processing times, then processing the tasks in decreasing (respectively, increasing) order of their processing times stochastically minimizes  wasted time due to re-processing provided we enforce the condition that the machine suffers not more than a single failure in the course of processing.  These two rules are commonly known as LPT and SPT, respectively, in the machine scheduling literature. The result just cited is interesting; the limitation on the number of machine failures permitted must be noted, but we argue in a later section that this model does capture some real-world situations. This result was extended to likelihood-ratio ordered processing times by \cite{righter94}.

The only paper in the  literature we have seen in which disruptions are modeled by an NHPP is \cite{pinedoross80}, who  consider a fixed number of tasks that must be attempted in some sequence. To successfully complete any given task requires a random amount of time, and  there are external  shocks, which occur according to an NHPP. If no shocks occur while a task is being performed, then that task is considered a success. If a shock does occur, then work on that task ends and work on the next one begins.  A fixed reward is obtained upon successful completion of a task. The authors determine the  sequence that  maximizes the expected number of successful tasks, the length of time until no tasks remain, and the expected total reward earned.
A stream of research that has some thematic overlap with the present work is to be found in queueing theoretic models of service systems in which the server is subject to random disruptions (for example, \cite{sengupta1990}, \cite{takine1997}). But whereas the central research problem in that stream of work is to characterize the distribution of waiting time - and other measures of service -  in a queueing system for a given  service protocol, our focus is on the determination of the optimal service discipline to minimize the makespan with a fixed set of jobs
available at time zero. The demand process, and the queueing behavior that it entails, lurk in the background of our model. 

\noindent {\bf Real-world settings:}  An NHPP with an increasing (decreasing) intensity function models a disruptive phase in which the rate of occurrence of disruptions is steadily increasing (decreasing) over time although the incidence of disruptions is random. It is reasonable to posit monotonically increasing disruption rates in the early phase of operation for machines that have infant mortality, or machines that have a burn-in period before stabilization of operations. On the other hand, machines that are in the last phase of their useful life may experience disruptions at a steadily increasing rate before requiring replacement. 
In \cite{bommerandfendley2018}, and \cite{marketal2002},  realizations of mistakes due to mental overload or underload can be modeled by an NHPP.  Mistakes in high mental overload situations can be modeled by an NHPP with an increasing failure intensity. Similarly, mistakes in mental underload situations could be modeled using decreasing failure rates. In data networks, the packet drop-rate over an end-to-end link increases during peak times. This increase can be modeled by an NHPP with increasing intensity. Because no message is recompiled by the receiving node until the last packet arrives, sequencing the packets so as to minimize message delivery time (i.e., makespan) will reduce delay. Simple sequencing rules, as we will discuss here, can be implemented on readily available hardware and software. 
We acknowledge that the assumption of monotonically increasing or decreasing disruption rates is likely to prevail only over particular phases in the  lifetime of a machine; in practice it is likely that disruption rates will fluctuate over time. We conducted several Monte Carlo simulation experiments to gauge the impact of consistently fluctuating disruption rates over the lifetime of a machine; the results are reported in the last section of the paper.

\noindent {\bf Technical challenges:}  There are several technical challenges in our model.
First, due to the non-homogeneous nature of NHPP's, the time needed does not
enjoy the memoryless property as in the case of standard Poisson disruptions.
The time needed to process the next task depends on how long it takes
to process all the previous tasks. The non-homogeneous nature creates
long-term memory, and it is in general the case in probability theory
that there are fewer available tools and theory for non-homogeneous stochastic processes.
Second, unlike the homogeneous case, the expected time processing tasks do not yield closed-form formulas
under the NHPP setting. As we will see later, they satisfy some integral equations, but analytical
solutions do not seem to be available. 
Third, in the general NHPP setting, there is no simple coupling method from probability theory
that one can use to show the optimal ordering of the tasks without quantifying the expected
processing times.

\noindent {\bf Summary of Results:}  
Our analysis and results make novel contributions.
Our  main result in the general model of failures is that for a sequence of short tasks, 
if the disruption rate is decreasing fast enough, then the optimal sequence is SPT
and if the disruption rate is increasing fast enough, then the optimal sequence is LPT. 
In the important case when there is exactly one (perhaps catastrophic) failure, we prove that the optimal sequence rides on the nature of the monotonicity in the density function of the time to the failure. If the said density function is strictly decreasing, then SPT is optimal; if it is strictly increasing, then LPT is optimal. 
Our analysis relies on deriving a chain of integral equations
that the expected processing times satisfy, and providing explicitly
computable upper and lower bounds for the expected processing times
by a careful analysis utilizing the integral equations and the properties
of NHPP. Our model contributes to the theory of sequencing tasks on machines subject to random disruptions. In a more general sense, our model and results suggest simple operational safeguards whereby the loss of systemic productivity engendered by processing failures can be minimized.
To the best of our knowledge, this is the first paper to study the problem when the timing of the disruptions is driven by a non-stationary point process. 

\section{Analysis and Results}\label{analysis}
We now give a precise mathematical description of our model, and then proceed to the analysis and results, beginning with the simplest case of two tasks. We will then extend the analysis to the general case. We relegate the proofs to the Appendix.

\subsection{Two Tasks}\label{sec:two}
In this section, we consider two tasks with lengths $a$ and $b$;
that is to say, if there were no disruptions, the time required to complete the tasks would be $a$ and $b$ time units, respectively. There is a non-homogeneous Poisson process with intensity $\lambda(t)$
that generates the random points in time at which disruptions occur. When a disruption occurs before a task is finished,
one has to go back to the beginning to restart the task. 

Let $\tau_{a,b}$ be the time to complete both tasks $a$ and $b$ given $a$ is processed before $b$. Define $\tau_{b,a}$, and for the single tasks, $\tau_{a}$ and $\tau_{b}$ similarly. Also define, assuming a disruption occurs at time $t$, the conditional mean remaining makespan: 
\begin{align*}
&M_{a,b}(t)=\mathbb{E}[\tau_{a,b}-t|\tau_{a}>t],
\\
&M_{b}(t)=\mathbb{E}[\tau_{b}-t|\tau_{b}>t\geq\tau_{a}],
\end{align*}
with $M_{b,a}(t)$ and $M_{a}(t)$ similarly defined. 
The following lemma, and similar equations for $M_{b,a}(t)$ and $M_{a}(t)$, follows easily from the properties of NHPP's, and the proof is omitted.

\begin{lemma}\label{lem:equation:two:tasks}
For every $t\geq 0$,
\begin{align}
M_{a,b}(t)&=(a+b)e^{-\int_{t}^{t+a+b}\lambda(s)ds}
+\int_{t}^{t+a}\lambda(s)e^{-\int_{t}^{s}\lambda(u)du}(s-t+M_{a,b}(s))ds
\nonumber
\\
&\qquad
+\int_{t+a}^{t+a+b}\lambda(s)e^{-\int_{t}^{s}\lambda(u)du}(s-t+M_{b}(s))ds,
\label{eqn:a:b}
\end{align}
where for every $t\geq 0$,
\begin{equation}\label{eqn:b}
M_{b}(t)=be^{-\int_{t}^{t+b}\lambda(s)ds}
+\int_{t}^{t+b}\lambda(s)e^{-\int_{t}^{s}\lambda(u)du}(s-t+M_{b}(s))ds.
\end{equation}
\end{lemma}


We are interested in comparing the values of $M_{a,b}(0)$ and $M_{b,a}(0)$.
There do not seem to be analytical solutions for $M_{a,b}(t)$ and $M_{b,a}(t)$
with explicit formulas except for the very special case $\lambda(t)\equiv\lambda$. 
In this case, none of the conditional mean makespans depend on $t$, so we denote them by $M_{a:\lambda}$ instead of $M_{a}(t)$, etc. 
Lemma~\ref{lem:const}, which will be used to provide bounds and control error terms for the more general NHPP case,
follows easily and the proof of Lemma~\ref{lem:const} is standard and thus omitted.

\begin{lemma}\label{lem:const}
(i) For all $a,\lambda>0$,
$M_{a;\lambda}=\frac{1}{\lambda}(e^{\lambda a}-1)$.

(ii) For all $a,b,\lambda>0$,
$M_{a,b;\lambda}=M_{b,a;\lambda}=\frac{1}{\lambda}(e^{\lambda a}+e^{\lambda b}-2)$.
\end{lemma}


Lemma~\ref{lem:const} will merely serve as a technical lemma that will be used
to provide bounds and control error terms later in the paper,
and in the rest of the paper, we will focus on the case when the disruption rates are non-constant functions of time. In this situation
the main challenge to compare $M_{a,b}(0)$ and $M_{b,a}(0)$ lies
in the fact that they do not have closed-form expressions.
However, in the next proposition, we show that under some special conditions,
one can nevertheless compute out $M_{a,b}(0)$ and $M_{b,a}(0)$ in closed-form
and hence make the comparison.


\begin{proposition}\label{prop:special}
Assume that $b>a>0$ and 
$\lambda(t)=0$ for $0\leq t\leq b$,
and $\lambda(t)=\lambda$ for $t>b$.
Then 
$M_{a,b}(0)=b+\frac{1}{\lambda}(e^{\lambda b}-e^{\lambda(b-a)})$,
$M_{b,a}(0)=b+\frac{1}{\lambda}(e^{\lambda a}-1)$,
and LPT is optimal.
\end{proposition}

Next, we will show that given a disruption rate $\lambda(t)$, $t\geq 0$,
when we have two tasks of length $a$ and $b$, where $a$ is much shorter
than $b$, then it is optimal to schedule $a$ first when $\lambda(t)$ is decreasing
and $b$ first when $\lambda(t)$ is increasing. 
Before we proceed, let us introduce some notation to facilitate the presentation and analysis. 
Let $\lambda(t)=\bar{\lambda} f(t)$, where $f(t)$ is a positive differentiable function
such that $f_{-}\leq f(t)\leq f_{+}$ for all $t$, where $f_{-},f_{+}$ are two positive constants,
and $\bar{\lambda}>0$ is a scaling parameter that denotes the scale of disruptions.

\begin{proposition}\label{prop:small:a}
Assume that $\lambda(t)$ is $L$-Lipschitz for some $L>0$.

(i) If $\lambda(t)$ is decreasing for all $t$ and strictly decreasing for $t\leq b$,
then SPT is optimal
given that
\begin{equation}
a<\min\left\{\frac{1}{2f_{+}\bar{\lambda}},(\mathcal{M}_{1})^{-1}
\left(1-e^{b\lambda(b)-\int_{0}^{b}\lambda(u)du}\right)\right\},
\end{equation}
where
\begin{align}
\mathcal{M}_{1}&:=
(L+(\bar{\lambda}f_{+})^{2})b
+(\bar{\lambda}f_{+})^{2}\left(b+\frac{e^{\lambda(0)b}-1}{\lambda(0)}\right)
\nonumber
\\
&\qquad
+\frac{3}{2}\lambda(0)\lambda(b)\left(\frac{1}{2\bar{\lambda}f_{+}}+b+\frac{e^{\lambda(b)b}-1}{\lambda(b)}\right)
+\lambda(b)+LM_{b}(0)+(\lambda(0))^{2}M_{b}(0)+\frac{3}{2}\lambda(0).\label{M1:defn}
\end{align}

(ii) If $\lambda(t)$ is increasing for all $t$ and strictly increasing for $t\leq b$,
then LPT is optimal
given that 
\begin{equation}
a<\min\left\{\frac{1}{2f_{+}\bar{\lambda}},(\mathcal{M}_{2})^{-1}e^{-\int_{0}^{b}\lambda(u)du}\left(e^{\lambda(b)b}-e^{\int_{0}^{b}\lambda(u)du}\right)\right\},
\end{equation}
where
\begin{align}\label{M2:defn}
\mathcal{M}_{2}:=\frac{5}{2}\bar{\lambda}f_{+}+bL
+\lambda(0)\left(e^{\bar{\lambda}f_{+}b}+\frac{1}{2}\right)
+\lambda(b)\left(b\bar{\lambda}f_{+}-1+e^{\bar{\lambda}f_{+}b}\right).
\end{align}
\end{proposition}


In Proposition~\ref{prop:small:a}, $a$ denotes the length of the shorter task. Thus 
Proposition~\ref{prop:small:a} infers that when you have two tasks of length $a$ and $b$, and $a$ is much shorter, 
then if $\lambda(t)$ is decreasing for every $t$ and strictly decreasing for $t\leq b$,
then SPT is optimal; on the other hand, if $\lambda(t)$ is increasing for every $t$ and strictly increasing for $t\leq b$,
then LPT is optimal.
The main challenge to prove Proposition~\ref{prop:small:a} is that $M_{a,b}(0)$ and $M_{b,a}(0)$
cannot be computed in closed-form. The proof instead relies on a first-order expansion in $a$ and some careful analysis to control all the error terms from the expansion.

\subsection{$n$ Tasks}\label{sec:n}

In this section, we extend our analysis to the case with $n$ tasks, of lengths $a_{1},a_{2},\ldots,a_{n}$,
and, without loss of generality, let us assume that $a_{1}<a_{2}<\cdots<a_{n}$.
Let $\tau_{1:n}$ be the time when all the tasks $a_{1},\ldots,a_{n}$ are finished
and the tasks are arranged in the order $a_{1},a_{2},\ldots,a_{n}$.
We define $\tau_{\pi(1):\pi(n)}$ similarly for the case when the tasks
are arranged in the order $a_{\pi(1)},a_{\pi(2)},\ldots,a_{\pi(n)}$,
where $\pi:\{1,\ldots,n\}\rightarrow\{1,\ldots,n\}$
is a permutation of $\{1,\ldots,n\}$.
Let us define:
\begin{align*}
&M_{1:n}(t)=\mathbb{E}\left[\tau_{1:n}-t|t<\tau_{a_{1}}\right],
\\
&M_{\pi(1):\pi(n)}(t)=\mathbb{E}\left[\tau_{\pi(1):\pi(n)}-t|t<\tau_{a_{\pi(1)}}\right].
\end{align*}
Similar to Lemma~\ref{lem:equation:two:tasks}, we have the following equation for $M_{1:n}(t)$, and the equation for $M_{\pi(1):\pi(n)}(t)$ is similar. For every $t\geq 0$,
\begin{align}
&M_{1:n}(t)=A_{n}e^{-\int_{t}^{t+A_{n}}\lambda(s)ds}
+\int_{t}^{t+A_{1}}\lambda(s)e^{-\int_{t}^{s}\lambda(u)du}(s-t+M_{1:n}(s))ds
\nonumber
\\
&\qquad\qquad\qquad\qquad\qquad\qquad
+\cdots+\int_{t+A_{n-1}}^{t+A_{n}}\lambda(s)e^{-\int_{t}^{s}\lambda(u)du}(s-t+M_{n:n}(s))ds,\label{eqn:a:1:n}
\end{align}
where $A_{i}:=a_{1}+a_{2}+\cdots+a_{i}$. 
To facilitate the presentation, we also define
$A_{i:j}:=\sum_{k=i}^{j}a_{k}$, 
$A_{\pi(i)}:=a_{\pi(1)}+a_{\pi(2)}+\cdots+a_{\pi(i)}$
and
$A_{\pi(i):\pi(j)}:=\sum_{k=i}^{j}a_{\pi(k)}$.
We are interested in comparing $M_{1:n}(0)$ with $M_{\pi(1):\pi(n)}(0)$.
Let us recall that
$\lambda(t)=\bar{\lambda} f(t)$, where $f(t)$ is a positive differentiable function such that $f_{-}\leq f(t)\leq f_{+}$
for every $t$, where $f_{-}=\inf_{t\geq 0}f(t)$, $f_{+}=\sup_{t\geq 0}f(t)$ are two positive constants,
and $\bar{\lambda}>0$ is a scaling parameter that denotes the scale of disruptions.
We will show that if $f(t)$ is decreasing,
then for any sufficiently small $\bar{\lambda}$, $M_{1:n}(0)<M_{\pi(1):\pi(n)}(0)$
for any permutation $\pi$ such that $(\pi(1),\ldots,\pi(n))\neq(1,2,\ldots,n)$,
and if $f(t)$ is increasing, 
then for any sufficiently small $\bar{\lambda}$, $M_{n:1}(0)<M_{\pi(1):\pi(n)}(0)$
for any permutation $\pi$ such that $(\pi(1),\ldots,\pi(n))\neq(n,n-1,\ldots,1)$.
Indeed, we can provide an explicitly computable upper bound
for $\bar{\lambda}$ such that the monotonicity results hold.

\begin{theorem}\label{thm:n}
(i) If $f(t)$ is strictly decreasing for $t\leq A_{n}$, then SPT is optimal
given $\bar{\lambda}\leq\frac{1}{2f_{+}a_{n}}$ and
\begin{equation}
\bar{\lambda}<\frac{
\sum_{i=1}^{n}a_{\pi(i)}\int_{0}^{A_{\pi(i)}}f(s)ds
-\sum_{i=1}^{n}a_{i}\int_{0}^{A_{i}}f(s)ds}
{(f_{+})^{2}A_{n}\sum_{i=1}^{n}a_{i}^{2}
+(f_{+})^{2}\frac{3}{4}\left(A_{n}\right)^{3}}.
\end{equation}

(ii) If $f(t)$ is strictly increasing for $t\leq A_{n}$, then 
LPT is optimal
given $\bar{\lambda}\leq\frac{1}{2f_{+}a_{n}}$ and
\begin{equation}
\bar{\lambda}<\frac{
\sum_{i=1}^{n}a_{\pi(i)}\int_{0}^{A_{\pi(i)}}f(s)ds
-\sum_{i=1}^{n}a_{n+1-i}\int_{0}^{A_{n+1-i:n}}f(s)ds}
{(f_{+})^{2}A_{n}\sum_{i=1}^{n}a_{i}^{2}
+(f_{+})^{2}\frac{3}{4}\left(A_{n}\right)^{3}}.
\end{equation}
\end{theorem}



\begin{remark}
There are situations in which not all sequences of tasks are feasible, owing to technological constraints. 
In these situations, Theorem~\ref{thm:n}---and indeed, all the results in this section---apply in the following modified sense.
Consider the subset $\Pi$ of processing times $a_{i_1}, a_{i_2},\cdots,a_{i_m}$ corresponding to the $m < n$ tasks that can be freely permuted amongst themselves, respecting the technological constraints. We say that a sequence of processing times is {\it relatively SPT(LPT)} when $\Pi$ is in SPT(LPT) order.  Theorem~\ref{thm:n} holds with SPT and LPT replaced by `relatively SPT' and `relatively LPT', respectively (noting that the main ingredient in the proof of Theorem~\ref{thm:n} is Lemma~\ref{lem:ineq:n} whose proof holds for $\Pi$, so that Theorem~\ref{thm:n} applies to this modified setup too). An example when this arises in actual practice is in project scheduling and assembly line processing, when substantial subsets of the sequential project network or assembly line follow a rigid sequence of tasks out of technological necessity.
\end{remark}

Our next result shows that we can relax the requirement that the disruption rate be sufficiently small.
Theorem~\ref{thm:n} holds for a sequence of short tasks provided that the disruption rate is increasing (or decreasing)
fast enough. Before we proceed, let us introduce the notation
$a_{i}:=\epsilon\bar{a}_{i}$, where $\bar{a}_{1}<\cdots<\bar{a}_{n}$ with $\bar{A}_{i}$, $\bar{A}_{i:j}$, 
$\bar{A}_{\pi(i)}$ and $\bar{A}_{\pi(i):\pi(j)}$ similarly defined as $A_{i},A_{i:j}, A_{\pi(i)}$ and $A_{\pi(i):\pi(j)}$.
We have the following result.

\begin{theorem}\label{thm:n:short}
(i) If $f(t)$ is strictly decreasing for $t\leq A_{n}$,
and for any permutation $\pi$ with $(\pi(1),\ldots,\pi(n))\neq(1,2,\ldots,n)$
\begin{equation}
|f'(0)|>\frac{2\bar{\lambda}\left((f_{+})^{2}\bar{A}_{n}\sum_{i=1}^{n}\bar{a}_{i}^{2}
+(f_{+})^{2}\frac{3}{4}\left(\bar{A}_{n}\right)^{3}\right)}{\sum_{i=1}^{n}\bar{a}_{i}\left(\bar{A}_{i}\right)^{2}-\sum_{i=1}^{n}\bar{a}_{\pi(i)}\left(\bar{A}_{\pi(i)}\right)^{2}},
\end{equation}
then SPT is optimal for sufficiently small $\epsilon>0$.

(ii) If $f(t)$ is strictly increasing for $t\leq A_{n}$,
and for any permutation $\pi$ with $(\pi(1),\ldots,\pi(n))\neq(n,n-1,\ldots,1)$
\begin{equation}
f'(0)>\frac{2\bar{\lambda}\left((f_{+})^{2}\bar{A}_{n}\sum_{i=1}^{n}\bar{a}_{i}^{2}
+(f_{+})^{2}\frac{3}{4}\left(\bar{A}_{n}\right)^{3}\right)}{\sum_{i=1}^{n}\bar{a}_{\pi(i)}\left(\bar{A}_{\pi(i)}\right)^{2}
-\sum_{i=1}^{n}\bar{a}_{n+1-i}\left(\bar{A}_{n+1-i:n}\right)^{2}},
\end{equation}
then LPT is optimal for sufficiently small $\epsilon>0$.
\end{theorem}


In Theorem~\ref{thm:n:short}, $M_{1:n}(0)$ and $M_{\pi(1):\pi(n)}(0)$
cannot be computed in closed-form, and the results of Theorem~\ref{thm:n:short}
are obtained by applying Theorem~\ref{thm:n} 
and a careful analysis for $\epsilon\rightarrow 0$ to check the conditions in Theorem~\ref{thm:n}.
Notice that in Theorem~\ref{thm:n}, we do not need $\lambda(t)$
to be strictly increasing (or decreasing) for all $t\geq 0$ and being strictly increasing (or decreasing)
for $0\leq t\leq A_{n}$ suffices. 
In general, if there exists some $t_{0}>0$
such that $\lambda(t)$ is strictly increasing (or decreasing)
for all $0\leq t\leq t_{0}$
and $\lambda(t)\equiv \lambda(t_{0})$ for all $t>t_{0}$, then the condition of being strictly increasing (or decreasing)
for $0\leq t\leq A_{n}$ in Theorem~\ref{thm:n}
holds if $A_{n}\leq t_{0}$, i.e. when the tasks are short. 
Next, we will show that when the tasks are long, 
the ordering will not matter, i.e. $M_{1:n}(0)=M_{\pi(1):\pi(n)}(0)$
for any permutation $\pi$ of $\{1,\ldots,n\}$.

\begin{theorem}\label{thm:equal:n}
Assume there exists some $t_{0}>0$
such that $\lambda(t)\equiv\lambda(t_{0})$ for all $t>t_{0}$.
If $a_{i}>t_{0}$ for every $1\leq i\leq n$, 
then $M_{1:n}(0)=M_{\pi(1):\pi(n)}(0)$
for any permutation $\pi$ of $\{1,\ldots,n\}$.
\end{theorem}

In other words, in the special case when the NHPP becomes a Poisson process with
constant rate $\lambda(t_{0})$ before any of the tasks
can complete, Theorem~\ref{thm:equal:n} shows that the order does not matter.
Theorem~\ref{thm:equal:n} is obtained by using the coupling method from probability theory
without knowing the exact formulas for $M_{1:n}(0)$
and $M_{\pi(1):\pi(n)}(0)$ which cannot be computed in closed-form.


\subsection{Single Failure Model}

In this section, we consider the same model, but with the constraint that
at most one failure can occur. We recall that this model was studied by \cite{Adiri1989}, with the difference
that they were interested in flowtime.

\begin{remark}
The one failure model captures real-world settings in which an error results in a one-time failure event---for instance, a catastrophe which is unlikely to be repeated. In such situations, it is imperative that the tasks be sequenced so as to minimize the expected loss of capacity. In this model, it is noteworthy that we are able to find concrete and easily identifiable conditions for the optimality of simple permutation schedules.
\end{remark}

Define $p(t):=\lambda(t)e^{-\int_{0}^{t}\lambda(u)du}$ as the probability density function of the first disruption time.
Note that if $\lambda(t)$ is differentiable, 
then $p'(t)=(\lambda'(t)-\lambda(t))e^{-\int_{0}^{t}\lambda(u)du}$
and thus $p(t)$ is strictly decreasing
if and only if $\lambda'(t)<\lambda(t)$
and $p(t)$ is strictly increasing 
if and only if $\lambda'(t)>\lambda(t)$.
Note that if $\lambda(t)$ is decreasing,
then we must have $\lambda'(t)<\lambda(t)$
and $p(t)$ is strictly decreasing.
But if $\lambda(t)$ is increasing, 
we still have $p(t)$ decreasing in $t$.
Therefore, in the original model (Section~\ref{sec:two} and Section~\ref{sec:n}), 
the criterion and cutoff is the monotonicity
of the disruption rate $\lambda(t)$,
whereas in the one breakdown model, 
the criterion and cutoff is the monotonicity
of the p.d.f. of the first disruption time $p(t)$.

Assume that there are $n$ tasks with lengths $a_{1},a_{2},\ldots,a_{n}$
and without loss of generality, let us assume that $a_{1}<a_{2}<\cdots<a_{n}$.
Let $T_{1:n}$ be the time when all the tasks $a_{1},\ldots,a_{n}$ are finished
with at most one breakdown
and the tasks are arranged in the order $a_{1},a_{2},\ldots,a_{n}$.
Similarly, we define $T_{\pi(1):\pi(n)}$ for the case when the tasks
with at most one breakdown
are arranged in the order $a_{\pi(1)},a_{\pi(2)},\ldots,a_{\pi(n)}$,
where $\pi:\{1,\ldots,n\}\rightarrow\{1,\ldots,n\}$
is a permutation of $\{1,\ldots,n\}$.
Let us define:
\begin{align*}
&R_{1:n}(t)=\mathbb{E}[T_{1:n}-t|\text{none of the tasks $a_{1},\ldots,a_{n}$ is finished at time $t$}],
\\
&R_{\pi(1):\pi(n)}(t)=\mathbb{E}[T_{\pi(1):\pi(n)}-t|\text{none of the tasks $a_{\pi(1)},\ldots,a_{\pi(n)}$ is finished at time $t$}].
\end{align*}
Similar as equation~\eqref{eqn:a:1:n}, we can write down the equation for $R_{1:n}(t)$, and the equation for $R_{\pi(1):\pi(n)}(t)$ is similar.
For every $t\geq 0$,
\begin{align}
R_{1:n}(t)
&=A_{n}e^{-\int_{t}^{t+A_{n}}\lambda(s)ds}
+\int_{t}^{t+A_{1}}\lambda(s)e^{-\int_{t}^{s}\lambda(u)du}(s-t+A_{n})ds
\nonumber
\\
&\quad
+\int_{t+A_{1}}^{t+A_{2}}\lambda(s)e^{-\int_{t}^{s}\lambda(u)du}(s-t+A_{2:n})ds
+\cdots+\int_{t+A_{n-1}}^{t+A_{n}}\lambda(s)e^{-\int_{t}^{s}\lambda(u)du}(s-t+A_{n:n})ds.
\label{R:a:1:n:t}
\end{align}
Finally, we characterize the difference of $R_{1:n}(0)$ and $R_{\pi(1):\pi(n)}(0)$ in a simple closed-form in the next result
and show that $R_{1:n}(0)$ is less than (greater than) $R_{\pi(1):\pi(n)}(0)$
depending on the monotonicity of the p.d.f. of the first disruption time $p(t)$.

\begin{theorem}\label{thm:one:break:n}
For any permutation $\pi$ of $\{1,2,\ldots,n\}$, we have
\begin{equation*}
R_{1:n}(0)
-R_{\pi(1):\pi(n)}(0)
=\sum_{i=1}^{n}\left(a_{i}\int_{0}^{A_{i}}p(s)ds-a_{\pi(i)}\int_{0}^{A_{\pi(i)}}p(s)ds\right).
\end{equation*}
Moreover, 
(i) If $p(t)$ is strictly decreasing for $t\leq A_{n}$, then
SPT is optimal; 
(ii) If $p(t)$ is strictly increasing for $t\leq A_{n}$, then 
LPT is optimal.
\end{theorem}

Unlike $M_{1:n}(0)$ in the previous section, we have a closed-form
expression for $R_{1:n}(0)$ in \eqref{R:a:1:n:t}, which makes the analysis 
for single failure model much simpler. Theorem~\ref{thm:one:break:n}
is obtained by using the expression \eqref{R:a:1:n:t} and some careful analysis.

\section{Experimental Results and Discussion}

We conducted Monte Carlo simulation experiments to gather empirical evidence about the performance of LPT and SPT for seven different types of rate functions $\lambda(t)$: linear increasing  ($\lambda(t) = \min(at, \lambda)$), concave increasing($\lambda(t)=\lambda \sqrt{at}$), increasing step function ($\lambda(t)= \lambda/2, t < t_0, ~\lambda ~\text{o/w} $), linear decreasing ($\lambda(t)= \max(\lambda - at,\lambda_0)$), convex decreasing ($\lambda(t)= 1/ \sqrt{t+1}$), decreasing step function ($\lambda(t)= \lambda, t < t_0, \lambda/2 ~\text{o/w} $), a non-monotonic function ($\lambda(t)=\lambda(1 + \sin{(at)})$), and a bathtub function composed of linear decreasing ($\lambda - at$), constant ($\lambda/2$), linear increasing ($\lambda/2 + at$) and constant functions ($\lambda$). We used the thinning of a Poisson process as well as \c{C}{\i}nlar's method \cite{Cinlar1975} to generate the NHPP failures.

For a batch of $n=4$ tasks with lengths $2, 4, 6$ and $8$ time units (so $A_{4}=20$), we simulated the disruption process for each permutation for $2,000,000$ replications (sample paths) and report the average makespan and its standard error in parentheses for $\lambda=0.4$ and appropriate values of $a$ and $t_0$ for monotone failure rates. We do not report the results for non-monotone rate functions, since the optimal sequence varies for different choices of $a$ and $t_0$. 




\begin{table}[htbp]
  \centering
  \begin{footnotesize}
    \begin{tabular}{c|c|c|c}
    \hline
 $\lambda(t)$  & \textbf{SPT batch time} &  \textbf{LPT batch time} & \textbf{Maximum $\%$ cost} \\
    & & & \textbf{of mis-sequencing}\\
    \hline
 Linear decreasing      & 41.34 (3.57E-3)  & 45.23 (7.71E-3) & 9.41 \\
 Step function decreasing    & 45.08 (1.14E-2)  & 48.54 (1.17E-2) & 7.68 \\
 Convex decreasing    & 27.67 (4.38E-3) & 29.01 (4.96E-3) & 4.84 \\
\hline
 Linear increasing       & 35.20 (2.17E-2)  & 26.50 (8.90 E-3) & 32.83 \\
 Step function increasing   & 86.88 (3.99E-2)   &  72.70 (3.63 E-2) & 19.50 \\
 Concave increasing     & 60.06 (3.46E-2) &  45.54 (2.49E-2) & 31.88  \\
    \hline
    \end{tabular}%
    \end{footnotesize}
     \caption{Experimental Results.}
  \label{table1}%
\end{table}%

Though our results in Section~\ref{analysis} for the optimality of SPT or LPT hold only for monotonone $\lambda(t)$ under restricted conditions, in our numerical results, we found that for all the monotone $\lambda(t)$ cases, SPT or LPT was optimal. In these cases, when SPT (LPT) is optimal, LPT (SPT) is the worst sequence. The difference in all cases is significant (note the standard errors in Table~\ref{table1}). On the other hand, neither LPT nor SPT were optimal for non-monotone rate functions we used. 
The cost of choosing the wrong sequence varied between $4.84$ percent to $32.83$ percent, depending on the form of the rate function (please see Table~\ref{table1}). That is to say, if one picks the SPT (LPT) sequence instead of the LPT (SPT) sequence in a situation in which the optimal sequence is LPT (SPT), then this sub-optimal choice resulted in a batch completion  time of between $4.84$ percent to $32.83$ percent more than the batch completion time associated with the optimal sequence. 

On the basis of the simulation results discussed above, we conjecture that when the processing of a fixed batch of tasks is disrupted by an NHPP with a monotone intensity function, SPT is optimal whenever the intensity function is decreasing and LPT is optimal whenever the intensity function is increasing. However, for the technical reasons that we gave at the start of Section~\ref{analysis}, proving this conjecture appears to be a very challenging problem.

{\bf A simple approximation to the original problem:} In order to give some concrete support for our  conjecture, we discuss a variation of our original problem with two tasks that approximates its main features in a much more tractable setup. The solution to this problem gives some intuitive  support for our conjecture about the connection between the monotonicity of the failure rate function and the optimality of SPT/LPT sequences. 

Assume that the first task experiences failures generated by a Poisson process with rate $\lambda_1$. Also assume that upon successful completion of the first task the failure process {\it instantaneously switches to a Poisson processes with rate $\lambda_2$} (this assumption allows us to compute the wasted time for each task independently).  Imposing the condition $\lambda_2 > \lambda_1$ this set up approximates an NHPP-like failure process with increasing failure rate.  

First we recall from Lemma~\ref{lem:const}(i) that the expected time $\mathbb{E}[\tau_t]$
to complete the task of length $t$ is
given by the explicit formula:
\begin{equation}\label{eqn:ew1}
\mathbb{E}[\tau_{t}]=\frac{e^{\lambda t}}{\lambda}-\frac{1}{\lambda}.
\end{equation}
Let $a < b$ be task lengths for two tasks. Then $\mathbb{E}[\tau_{a,b}]$ and $\mathbb{E}[\tau_{b,a}]$ are the expected completion times of the SPT and LPT sequences, respectively.
Also, let $\lambda_1 < \lambda_2$. \emph{Given this setup, we expect that LPT is the optimal two-task sequence}. Using Equation~\eqref{eqn:ew1} we can write
\begin{align*}
\mathbb{E}[\tau_{a,b}] = \frac{e^{\lambda_1 a}}{\lambda_1} - \frac{1}{\lambda_1} + \frac{e^{\lambda_2 b}}{\lambda_2} - \frac{1}{\lambda_2}, 
\qquad
\mathbb{E}[\tau_{b,a}] = \frac{e^{\lambda_1 b}}{\lambda_1}- \frac{1}{\lambda_1} + \frac{e^{\lambda_2 a}}{\lambda_2} - \frac{1}{\lambda_2},
\end{align*}
such that
\begin{align*}
\mathbb{E}[\tau_{a,b}] - \mathbb{E}[\tau_{b,a}]  &= \frac{e^{\lambda_1 a} - e^{\lambda_1 b }}{\lambda_1} -  \frac{e^{\lambda_2 a} - e^{\lambda_2 b} }{\lambda_2}. 
\end{align*}
Let $\lambda_2 = \lambda_1 + \delta$ with $\delta>0$ and $b=a+\epsilon$ with $\epsilon > 0$.
We have
\begin{align*}
\Delta =\mathbb{E}[\tau_{a,b}] - \mathbb{E}[\tau_{b,a}] 
=\frac{e^{\lambda_1 a} - e^{\lambda_1 a }e^{\lambda_1 \epsilon }}{\lambda_1} -  \frac{e^{\lambda_1 a}e^{\delta a} - e^{\lambda_1 a}e^{\delta a}e^{\lambda_1 \epsilon}e^{\delta \epsilon} }{\lambda_1 + \delta}.
\end{align*}
We first note that $\Delta=0$ when $\delta=0, \epsilon=0$. We now show that $\Delta$ is increasing in $\epsilon$ for every value of $\delta$ by computing that:
\begin{align}
\frac{d \Delta}{d\epsilon} = e^{\lambda_1 a}\left[ - \frac{\lambda_1 e^{\lambda_1 \epsilon }}{\lambda_1} + \frac{e^{\delta a} (\lambda_1 + \delta) e^{(\lambda_1+\delta) \epsilon}}{\lambda_1 + \delta}  \right] = e^{\lambda_1 a}\left[ -  e^{\lambda_1 \epsilon } + e^{\delta a} e^{(\lambda_1+\delta) \epsilon}  \right] \geq 0, \quad\text{for every $\delta,\epsilon \geq 0$.}\label{eq:dde}
\end{align}
Inequality \eqref{eq:dde} is trivially satisfied, which proves the result.

We have shown that LPT is optimal for this two-task scenario. Likewise, if $\lambda_2 < \lambda_1$ we can show that SPT is optimal. This illustrates that processing the longer (shorter) task when the failure rates are low (high) is better.

\section{Conclusions}

In this paper, we study the optimal sequencing
of a batch of tasks on a machine subject to 
random NHPP disruptions, where the intensity of NHPP
is monotonic in time. We find conditions
under which either SPT or LPT minimizes the expected
time that is needed to complete a batch of tasks.
We conjecture that it is---without any restrictions---optimal to apply either SPT or LTP sequencing rules \emph{when the NHPP intensity function is monotonic}. 
In our model, we imposed relatively strong
structure, and did not for example consider deterministic or i.i.d. repair times,
which would make the analysis even more technically
challenging. We leave this for future research.

\section*{Acknowledgements}

The authors thank Prof. Jayaram Sethuraman for helpful discussions. Lingjiong Zhu acknowledges the support from NSF Grants DMS-1613164, 
DMS-2053454 and a Simons Foundation Collaboration Grant.


\bibliographystyle{alpha} 
\bibliography{research} 

\newcommand{\etalchar}[1]{$^{#1}$}
\begin{thebibliography}{BFMK90}

\bibitem[ABFK89]{Adiri1989}
Igal Adiri, John Bruno, Esther Frostig, and A.H.G.~Rinnooy Kan.
\newblock Single machine flow-time scheduling with a single breakdown.
\newblock {\em Acta Informatica}, 26:679--696, 1989.
\newblock 

\bibitem[AFL{\etalchar{+}}08]{asmussenetal2008}
Soeren Asmussen, Pierre Fiorini, Lester Lipsky, Tomasz Rolski, and Robert
  Sheahan.
\newblock Asymptotic behavior of total times for jobs that must start over if a
  failure occurs.
\newblock {\em Mathematics of Operations Research}, 33(4):932--944, 2008.
\newblock 

\bibitem[AP12]{aytugandpaul2012}
Haldun Aytug and Anand Paul.
\newblock Sequencing jobs on a non-{M}arkovian machine with random disruptions.
\newblock {\em IIE Transactions}, 44(8):671--680, 2012.

\bibitem[BF18]{bommerandfendley2018}
Sharon~Claxton Bommer and Mary Fendley.
\newblock A theoretical framework for evaluating mental workload resources in
  human systems design for manufacturing operations.
\newblock {\em International Journal of Industrial Ergonomics}, 63:7--17, 2018.

\bibitem[BFMK90]{Birge1990}
J.~Birge, J.B.G. Frenk, J.~Mittenthal, and A.H.G.~Rinnooy Kan.
\newblock Single-machine scheduling subject to stochastic breakdowns.
\newblock {\em Naval Research Logistics}, 37:661--677, 1990.
\newblock 

\bibitem[\c{C}75]{Cinlar1975}
Erhan \c{C}{\i}nlar.
\newblock {\em Introduction to Stochastic Processes}.
\newblock Prentice Hall, New Jersey, 1975.

\bibitem[CSZ03]{Cai2003}
Xiaoqiang Cai, Xiaoqian Sun, and Xian Zhou.
\newblock Stochastic scheduling with preemptive-repeat machine breakdowns to
  minimize the expected weighted flow time.
\newblock {\em Probability in the Engineering and Informational Sciences},
  17(04):467--485, 2003.

\bibitem[CSZ04]{Cai2004}
Xiaoqiang Cai, Xiaoqian Sun, and Xian Zhou.
\newblock Stochastic scheduling subject to machine breakdowns: The
  preemptive-repeat model with discounted reward and other criteria.
\newblock {\em Naval Research Logistics}, 51:800--817, 2004.

\bibitem[Fro91]{Frostig1991}
Esther Frostig.
\newblock A note on stochastic scheduling on a single machine subject to
  breakdown - the preemptive repeat model.
\newblock {\em Probability in the Engineering and Informational Sciences},
  5:349--354, 1991.
\newblock 

\bibitem[KAP06]{Kasap2006}
Nihat Kasap, Haldun Aytug, and Anand Paul.
\newblock Minimizing makespan on a single machine subject to random breakdowns.
\newblock {\em Operations Research Letters}, 34(1):29--36, 2006.
\newblock 

\bibitem[KAP08]{Kasap2008}
Nihat Kasap, Haldun Aytug, and Anand Paul.
\newblock Erratum to ``{M}inimizing makespan on a single machine subject to
  random breakdowns''.
\newblock {\em Operations Research Letters}, 36:140, 2008.

\bibitem[PR80]{pinedoross80}
Michael~L. Pinedo and Sheldon~M. Ross.
\newblock Scheduling jobs subject to nonhomogeneous {P}oisson shocks.
\newblock {\em Management Science}, 26(12):1250--1257, 1980.

\bibitem[Rig94]{righter94}
Rhonda Righter.
\newblock Scheduling.
\newblock In Moshe Shaked and J.~George Shanthikumar, editors, {\em Stochastic
  Orders and Their Applications}, pages 381--432. Academic Press, 1994.

\bibitem[Sen90]{sengupta1990}
Bhaskar Sengupta.
\newblock A queue with service interruptions in an alternating random
  environment.
\newblock {\em Operations Research}, 38(2):308--318, 1990.
\newblock 

\bibitem[TS97]{takine1997}
Tetsuya Takine and Bhaskar Sengupta.
\newblock A single server queue with service interruptions.
\newblock {\em Queueing Systems}, 26:285--300, 1997.
\newblock 

\bibitem[YS02]{marketal2002}
Mark~S. Young and Neville~A. Stanton.
\newblock Attention and automation: New perspectives on mental underload and
  performance.
\newblock {\em Theoretical Issues in Ergonomics Science}, 3(2):178--194, 2002.

\end{thebibliography}







%




\newpage
\appendix

\section{Technical Proofs}\label{appendix:main}


\subsection{Proofs of Proposition~\ref{prop:special} and Proposition~\ref{prop:small:a}}

\noindent
{\bf Proof of Proposition~\ref{prop:special}}
Let us first compute $M_{b,a}$.
In this case, since $\lambda(t)=0$ for $t\leq b$, the time
to process $b$ is precisely $b$. After this task is finished,
the intensity is $\lambda(t)=\lambda$, and the expected value of 
the time to finish task $a$ is given by $\frac{1}{\lambda}(e^{\lambda a}-1)$.
Hence, we have
\begin{equation*}
M_{b,a}(0)=b+\frac{1}{\lambda}(e^{\lambda a}-1).
\end{equation*}
Next, let us compute $M_{a,b}(0)$.
Since $\lambda(t)=0$ for $t\leq b$ and we assumed that $a<b$, thus
it takes precisely time $a$ to finish the task $a$. 
Now, we need to finish task $b$. Between time $t=a$ and $t=b$, we
still have $\lambda(t)=0$. Therefore at least the amount $b-a$ out of task $b$
can be finished without disruption. The time to finish task $b$ is precisely $b$
when in the remaining $b-(b-a)=a$ length of time, there is no disruption, which happens
with probability $e^{-a\lambda}$. On the other hand, between time $t=b-a$ and $t=b$, 
there can be a Poisson disruption with intensity $\lambda$, and once it gets disrupted,
the expected time to finish the task would be $\frac{1}{\lambda}(e^{\lambda b}-1)$.
That is,
\begin{equation*}
M_{a,b}(0)=a+be^{-a\lambda}+\int_{0}^{a}\left(b-a+t+\frac{1}{\lambda}(e^{\lambda b}-1)\right)\lambda e^{-\lambda t}dt=b+\frac{1}{\lambda}(e^{\lambda b}-e^{\lambda(b-a)}).
\end{equation*}
Hence,
\begin{equation*}
M_{b,a}(0)-M_{a,b}(0)=\frac{1}{\lambda}(e^{\lambda a}-1-e^{\lambda b}+e^{\lambda(b-a)}).
\end{equation*}
Let us define
\begin{equation*}
G(b):=\frac{1}{\lambda}(e^{\lambda a}-1-e^{\lambda b}+e^{\lambda(b-a)}),
\qquad b\geq a.
\end{equation*}
It is easy to see that $G(a)=0$ and $G(b)$ is strictly decreasing in $b$,
which implies that $G(b)<0$ for any $b>a$.
Hence, we conclude that for $b>a$, and for this particular choice of $\lambda(t)$, we have
$M_{b,a}(0)<M_{a,b}(0)$.
\hfill $\Box$


\noindent
{\bf Proof of Proposition~\ref{prop:small:a}}
(i) For any $s\geq 0$ and $a\leq\frac{1}{2f_{+}\bar{\lambda}}$,
\begin{equation*}
a\leq M_{a}(s)\leq M_{a;\bar{\lambda}f_{+}}
=\frac{1}{\bar{\lambda}f_{+}}(e^{\bar{\lambda}f_{+}a}-1)
\leq a+\bar{\lambda}f_{+}a^{2},
\end{equation*}
and
\begin{equation*}
a+M_{b}(s)\leq M_{b,a}(s)\leq a+\bar{\lambda}f_{+}a^{2}+M_{b}(s).
\end{equation*}
The term $M_{a,b}(0)$ is the sum
of $M_{a}(0)$ and then the expected time
to complete task $b$ after finishing task $a$. 
When $\lambda(t)$ is decreasing, 
we have
\begin{equation*}
M_{a,b}(0)\leq M_{a}(0)+M_{b}(a)\leq a+\bar{\lambda}f_{+}a^{2}+M_{b}(a).
\end{equation*}
Hence, we get
\begin{equation}\label{by:connection}
M_{a,b}(0)-M_{b,a}(0)
\leq
\bar{\lambda}f_{+}a^{2}+M_{b}(a)-M_{b}(0).
\end{equation}

On the other hand, we can compute that
\begin{align*}
M_{b}(a)-M_{b}(0)
&=be^{-\int_{a}^{a+b}\lambda(s)ds}
-be^{-\int_{0}^{b}\lambda(s)ds}
\\
&\qquad
+\int_{a}^{a+b}\lambda(s)e^{-\int_{a}^{s}\lambda(u)du}(s-a+M_{b}(s))ds
-\int_{0}^{b}\lambda(s)e^{-\int_{0}^{s}\lambda(u)du}(s+M_{b}(s))ds
\\
&=be^{-\int_{0}^{b}\lambda(s)ds}\left(e^{-\int_{b}^{a+b}\lambda(s)ds+\int_{0}^{a}\lambda(s)ds}-1\right)
\\
&\quad
-a\int_{a}^{a+b}\lambda(s)e^{-\int_{a}^{s}\lambda(u)du}ds
\nonumber
\\
&\quad\quad
+\int_{a}^{a+b}\lambda(s)\left(e^{\int_{0}^{a}\lambda(u)du}-1\right)e^{-\int_{0}^{s}\lambda(u)du}(s+M_{b}(s))ds
\\
&\quad\quad\quad
+\int_{b}^{a+b}\lambda(s)e^{-\int_{0}^{s}\lambda(u)du}(s+M_{b}(s))ds
-\int_{0}^{a}\lambda(s)e^{-\int_{0}^{s}\lambda(u)du}(s+M_{b}(s))ds.
\end{align*}
For any $0\leq x\leq\frac{1}{2}$, by Taylor's expansion,
\begin{equation}\label{by:Taylor}
e^{x}-1=\sum_{k=1}^{\infty}\frac{x^{k}}{k!}
=x+\sum_{k=2}^{\infty}\frac{x^{k}}{k!}
\leq
x+\frac{1}{2}\sum_{k=2}^{\infty}x^{k}
=x+\frac{x^{2}}{2(1-x)}
\leq
x+x^{2}.
\end{equation}
Since $\lambda(t)$ is decreasing, 
$-\int_{b}^{a+b}\lambda(s)ds+\int_{0}^{a}\lambda(s)ds\geq 0$, 
and for any $a\leq\frac{1}{2\bar{\lambda}f_{+}}$, we have
\begin{equation}\label{by:one:half}
-\int_{b}^{a+b}\lambda(s)ds+\int_{0}^{a}\lambda(s)ds
\leq
\int_{0}^{a}\lambda(s)ds
\leq a\bar{\lambda}f_{+}\leq\frac{1}{2},
\end{equation}
so that by applying \eqref{by:Taylor} and \eqref{by:one:half} it follows that
\begin{align*}
&be^{-\int_{0}^{b}\lambda(s)ds}\left(e^{-\int_{b}^{a+b}\lambda(s)ds+\int_{0}^{a}\lambda(s)ds}-1\right)
\\
&\leq
be^{-\int_{0}^{b}\lambda(s)ds}
\left(-\int_{b}^{a+b}\lambda(s)ds+\int_{0}^{a}\lambda(s)ds
+\left(-\int_{b}^{a+b}\lambda(s)ds+\int_{0}^{a}\lambda(s)ds\right)^{2}\right)
\\
&\leq
be^{-\int_{0}^{b}\lambda(s)ds}
\left(-\int_{b}^{a+b}\lambda(s)ds+\int_{0}^{a}\lambda(s)ds
+(\bar{\lambda}f_{+})^{2}a^{2}\right)
\\
&\leq
be^{-\int_{0}^{b}\lambda(s)ds}
\left(-a\lambda(a+b)+a\lambda(0)+(\bar{\lambda}f_{+})^{2}a^{2}\right)
\\
&\leq
be^{-\int_{0}^{b}\lambda(s)ds}
\left(-a\lambda(b)+a\lambda(0)+La^{2}+(\bar{\lambda}f_{+})^{2}a^{2}\right)
\\
&=a(\lambda(0)-\lambda(b))be^{-\int_{0}^{b}\lambda(s)ds}
+a^{2}(L+(\bar{\lambda}f_{+})^{2})be^{-\int_{0}^{b}\lambda(s)ds}.
\end{align*}
Moreover, since $\lambda(s)$ is decreasing,
we have $\lambda(s)e^{-\int_{0}^{s}\lambda(u)du}$ is decreasing and  
\begin{align*}
-a\int_{a}^{a+b}\lambda(s)e^{-\int_{a}^{s}\lambda(u)du}ds
&\leq
-a\int_{a}^{a+b}\lambda(s)e^{-\int_{0}^{s}\lambda(u)du}ds
\\
&\leq
-a\int_{0}^{b}\lambda(s)e^{-\int_{0}^{s}\lambda(u)du}ds.
\end{align*}
Next, for any $a\leq\frac{1}{2\bar{\lambda}f_{+}}$, 
we have $\int_{0}^{a}\lambda(u)du\leq\frac{1}{2}$ so that
\begin{align}
&\int_{a}^{a+b}\lambda(s)\left(e^{\int_{0}^{a}\lambda(u)du}-1\right)e^{-\int_{0}^{s}\lambda(u)du}(s+M_{b}(s))ds
\nonumber
\\
&\leq
\left(a\lambda(0)+(\bar{\lambda}f_{+})^{2}a^{2}\right)\int_{a}^{a+b}\lambda(s)e^{-\int_{0}^{s}\lambda(u)du}(s+M_{b}(s))ds
\nonumber
\\
&
\leq
\left(a\lambda(0)+(\bar{\lambda}f_{+})^{2}a^{2}\right)\left(\int_{0}^{b}\lambda(s)e^{-\int_{0}^{s}\lambda(u)du}(s+M_{b}(s))ds
+\int_{b}^{a+b}\lambda(s)e^{-\int_{0}^{s}\lambda(u)du}(s+M_{b}(s))ds\right)
\nonumber
\\
&\leq
a\lambda(0)\int_{0}^{b}\lambda(s)e^{-\int_{0}^{s}\lambda(u)du}(s+M_{b}(s))ds
+(\bar{\lambda}f_{+})^{2}a^{2}(b+M_{b}(0))
\nonumber
\\
&\qquad
+\left(a\lambda(0)+(\bar{\lambda}f_{+})^{2}a^{2}\right)a\lambda(b)(a+b+M_{b}(b))
\label{explain:1}
\\
&\leq
a\lambda(0)\int_{0}^{b}\lambda(s)e^{-\int_{0}^{s}\lambda(u)du}(s+M_{b}(s))ds
+(\bar{\lambda}f_{+})^{2}a^{2}\left(b+\frac{e^{\lambda(0)b}-1}{\lambda(0)}\right)
\nonumber
\\
&\qquad\qquad\qquad
+\frac{3}{2}a^{2}\lambda(0)\lambda(b)\left(\frac{1}{2\bar{\lambda}f_{+}}+b+\frac{e^{\lambda(b)b}-1}{\lambda(b)}\right),
\label{explain:2}
\end{align}
where we used the fact that 
$M_{b}(s)$ is decreasing in $s$ (since $\lambda(s)$ is decreasing)
so that
for any $0\leq s\leq b$, $s+M_{b}(s)\leq b+M_{b}(0)$
and for any $b\leq s\leq a+b$, $s+M_{b}(s)\leq a+b+M_{b}(b)$
to obtain \eqref{explain:1}, 
and the fact that $M_{b}(0)\leq\frac{e^{\lambda(0)b}-1}{\lambda(0)}$ (since $\lambda(s)\leq\lambda(0)$
for any $s\geq 0$ and we can then apply Lemma~\ref{lem:const}) 
and the fact that $M_{b}(b)\leq\frac{e^{\lambda(b)b}-1}{\lambda(b)}$ (since $\lambda(s)\leq\lambda(b)$
for any $s\geq b$ and we can then apply Lemma~\ref{lem:const}) 
and finally the fact that $a\leq\frac{1}{2\bar{\lambda}f_{+}}$
and $\lambda(0)=\bar{\lambda}f_{+}$ 
to obtain \eqref{explain:2}.

Moreover, since $\lambda(t)$ is decreasing, we have
\begin{align*}
\int_{b}^{a+b}\lambda(s)e^{-\int_{0}^{s}\lambda(u)du}(s+M_{b}(s))ds
&\leq
a\lambda(b)e^{-\int_{0}^{b}\lambda(u)du}(a+b+M_{b}(b))
\\
&\leq
a\lambda(b)e^{-\int_{0}^{b}\lambda(u)du}(b+M_{b}(b))+a^{2}\lambda(b).
\end{align*}
Finally, 
\begin{align*}
&-\int_{0}^{a}\lambda(s)e^{-\int_{0}^{s}\lambda(u)du}(s+M_{b}(s))ds
\\
&\leq
-e^{-\int_{0}^{a}\lambda(u)du}\int_{0}^{a}\lambda(s)M_{b}(s)ds
\\
&\leq
-\left(1-\int_{0}^{a}\lambda(u)du\right)a\lambda(a)M_{b}(a)
\\
&\leq
-a\lambda(a)M_{b}(a)+a^{2}(\lambda(0))^{2}M_{b}(0)
\\
&=-a\lambda(a)M_{b}(0)+a\lambda(a)(M_{b}(0)-M_{b}(a))+a^{2}(\lambda(0))^{2}M_{b}(0)
\\
&\leq -a\lambda(0)M_{b}(0)+a^{2}LM_{b}(0)+a\lambda(0)(M_{b}(0)-M_{b}(a))+a^{2}(\lambda(0))^{2}M_{b}(0).
\end{align*}
Putting everything together, we get
\begin{align*}
&M_{b}(a)-M_{b}(0)
\\
&\leq
a(\lambda(0)-\lambda(b))be^{-\int_{0}^{b}\lambda(s)ds}
+a^{2}(L+(\bar{\lambda}f_{+})^{2})be^{-\int_{0}^{b}\lambda(s)ds}
-a\int_{0}^{b}\lambda(s)e^{-\int_{0}^{s}\lambda(u)du}ds
\\
&\quad
+a\lambda(0)\int_{0}^{b}\lambda(s)e^{-\int_{0}^{s}\lambda(u)du}(s+M_{b}(s))ds
+(\bar{\lambda}f_{+})^{2}a^{2}\left(b+\frac{e^{\lambda(0)b}-1}{\lambda(0)}\right)
\\
&\qquad
+\frac{3}{2}a^{2}\lambda(0)\lambda(b)\left(\frac{1}{2\bar{\lambda}f_{+}}+b+\frac{e^{\lambda(b)b}-1}{\lambda(b)}\right)
+a\lambda(b)e^{-\int_{0}^{b}\lambda(u)du}(b+M_{b}(b))+a^{2}\lambda(b)
\\
&
\qquad\qquad
-a\lambda(0)M_{b}(0)+a^{2}LM_{b}(0)+a\lambda(0)(M_{b}(0)-M_{b}(a))+a^{2}(\lambda(0))^{2}M_{b}(0).
\end{align*}
Let us recall that $a\lambda(0)=a\bar{\lambda}f_{+}\leq\frac{1}{2}$ and
\begin{equation*}
M_{b}(0)=be^{-\int_{0}^{b}\lambda(s)ds}
+\int_{0}^{b}\lambda(s)e^{-\int_{0}^{s}\lambda(u)du}(s+M_{b}(s))ds.
\end{equation*}
Therefore, we get
\begin{align*}
&\frac{3}{2}\left(M_{b}(a)-M_{b}(0)\right)
\\
&\leq
a\lambda(b)e^{-\int_{0}^{b}\lambda(u)du}M_{b}(b)-a\int_{0}^{b}\lambda(s)e^{-\int_{0}^{s}\lambda(u)du}ds
\\
&\quad
+a^{2}(L+(\bar{\lambda}f_{+})^{2})be^{-\int_{0}^{b}\lambda(s)ds}
+(\bar{\lambda}f_{+})^{2}a^{2}\left(b+\frac{e^{\lambda(0)b}-1}{\lambda(0)}\right)
\\
&\qquad
+\frac{3}{2}a^{2}\lambda(0)\lambda(b)\left(\frac{1}{2\bar{\lambda}f_{+}}+b+\frac{e^{\lambda(b)b}-1}{\lambda(b)}\right)
+a^{2}\lambda(b)
+a^{2}LM_{b}(0)+a^{2}(\lambda(0))^{2}M_{b}(0).
\end{align*}
Since $\lambda(t)$ is strictly decreasing for $0\leq t\leq b$, we can compute that
\begin{align*}
a\lambda(b)e^{-\int_{0}^{b}\lambda(u)du}M_{b}(b)-a\int_{0}^{b}\lambda(s)e^{-\int_{0}^{s}\lambda(u)du}ds
&=ae^{-\int_{0}^{b}\lambda(u)du}
\left(\lambda(b)M_{b}(b)+1-e^{\int_{0}^{b}\lambda(u)du}\right)
\\
&\leq
ae^{-\int_{0}^{b}\lambda(u)du}
\left(e^{b\lambda(b)}-e^{\int_{0}^{b}\lambda(u)du}\right)<0.
\end{align*}
By \eqref{by:connection}, we have
\begin{align*}
&M_{a,b}(0)-M_{b,a}(0)
\\
&\leq M_{b}(a)-M_{b}(0)+\bar{\lambda}f_{+}a^{2}
\\
&\leq
\frac{2}{3}ae^{-\int_{0}^{b}\lambda(u)du}
\left(e^{b\lambda(b)}-e^{\int_{0}^{b}\lambda(u)du}\right)
+\frac{2}{3}a^{2}(L+(\bar{\lambda}f_{+})^{2})b
+\frac{2}{3}(\bar{\lambda}f_{+})^{2}a^{2}\left(b+\frac{e^{\lambda(0)b}-1}{\lambda(0)}\right)
\\
&\qquad
+a^{2}\lambda(0)\lambda(b)\left(\frac{1}{2\bar{\lambda}f_{+}}+b+\frac{e^{\lambda(b)b}-1}{\lambda(b)}\right)
\\
&\qquad
+\frac{2}{3}a^{2}\lambda(b)
+\frac{2}{3}a^{2}LM_{b}(0)+\frac{2}{3}a^{2}(\lambda(0))^{2}M_{b}(0)+\lambda(0)a^{2}<0,
\end{align*}
provided that
\begin{equation*}
a<(\mathcal{M}_{1})^{-1}e^{-\int_{0}^{b}\lambda(u)du}
\left(e^{\int_{0}^{b}\lambda(u)du}-e^{b\lambda(b)}\right),
\end{equation*}
where $\mathcal{M}_{1}$ is defined in \eqref{M1:defn}.
This completes the proof of (i).

(ii) For any $s\geq 0$ and $a\leq\frac{1}{2f_{+}\bar{\lambda}}$,
\begin{equation*}
a\leq M_{a}(s)\leq M_{a;\bar{\lambda}f_{+}}
=\frac{1}{\bar{\lambda}f_{+}}(e^{\bar{\lambda}f_{+}a}-1)
\leq a+\bar{\lambda}f_{+}a^{2},
\end{equation*}
and
\begin{equation*}
a+M_{b}(s)\leq M_{b,a}(s)\leq a+\bar{\lambda}f_{+}a^{2}+M_{b}(s).
\end{equation*}
The term $M_{a,b}(0)$ is the sum
of $M_{a}(0)$ and then the expected time
to complete task $b$ after finishing task $a$. 
When $\lambda(t)$ is increasing, 
we have
$M_{a,b}(0)\geq M_{a}(0)+M_{b}(a)\geq a+M_{b}(a)$.
Hence, we get
\begin{equation*}
M_{a,b}(0)-M_{b,a}(0)
\geq
M_{b}(a)-M_{b}(0)-\bar{\lambda}f_{+}a^{2}.
\end{equation*}
It follows from the proof of (i) that
\begin{align*}
&M_{b}(a)-M_{b}(0)
\\
&=be^{-\int_{0}^{b}\lambda(s)ds}\left(e^{-\int_{b}^{a+b}\lambda(s)ds+\int_{0}^{a}\lambda(s)ds}-1\right)
\\
&\qquad
-a\int_{a}^{a+b}\lambda(s)e^{-\int_{a}^{s}\lambda(u)du}ds
+\int_{a}^{a+b}\lambda(s)\left(e^{\int_{0}^{a}\lambda(u)du}-1\right)e^{-\int_{0}^{s}\lambda(u)du}(s+M_{b}(s))ds
\\
&\qquad
+\int_{b}^{a+b}\lambda(s)e^{-\int_{0}^{s}\lambda(u)du}(s+M_{b}(s))ds
-\int_{0}^{a}\lambda(s)e^{-\int_{0}^{s}\lambda(u)du}(s+M_{b}(s))ds.
\end{align*}
Since $\lambda(t)$ is increasing, $\int_{b}^{a+b}\lambda(s)ds\geq\int_{0}^{a}\lambda(s)ds$, 
and
\begin{align*}
be^{-\int_{0}^{b}\lambda(s)ds}\left(e^{-\int_{b}^{a+b}\lambda(s)ds+\int_{0}^{a}\lambda(s)ds}-1\right)
&\geq
be^{-\int_{0}^{b}\lambda(s)ds}\left(\int_{b}^{a+b}\lambda(s)ds-\int_{0}^{a}\lambda(s)ds\right)
\\
&\geq
be^{-\int_{0}^{b}\lambda(s)ds}\left(a\lambda(b)-a\lambda(a)\right)
\\
&\geq
be^{-\int_{0}^{b}\lambda(s)ds}\left(a\lambda(b)-a\lambda(0)\right)
-be^{-\int_{0}^{b}\lambda(s)ds}a^{2}L.
\end{align*}
Next, we have
\begin{align*}
-a\int_{a}^{a+b}\lambda(s)e^{-\int_{a}^{s}\lambda(u)du}ds
&=-a\left(1-e^{-\int_{a}^{a+b}\lambda(u)du}\right)
\\
&=-a\left(1-e^{-\int_{0}^{b}\lambda(u)du}\right)
+a\left(e^{-\int_{a}^{a+b}\lambda(u)du}-e^{-\int_{0}^{b}\lambda(u)du}\right)
\\
&=-a\left(1-e^{-\int_{0}^{b}\lambda(u)du}\right)
+ae^{-\int_{0}^{b}\lambda(u)du}\left(e^{-\int_{b}^{a+b}\lambda(u)du+\int_{0}^{a}\lambda(u)du}-1\right)
\\
&\geq
-a\left(1-e^{-\int_{0}^{b}\lambda(u)du}\right)
-ae^{-\int_{0}^{b}\lambda(u)du}\left(\int_{a}^{a+b}\lambda(u)du-\int_{0}^{a}\lambda(u)du\right)
\\
&\geq
-a\left(1-e^{-\int_{0}^{b}\lambda(u)du}\right)
-a^{2}\bar{\lambda}f_{+}.
\end{align*}
Moreover, we have
\begin{align*}
&\int_{a}^{a+b}\lambda(s)\left(e^{\int_{0}^{a}\lambda(u)du}-1\right)e^{-\int_{0}^{s}\lambda(u)du}(s+M_{b}(s))ds
\\
&\geq
a\lambda(0)\int_{a}^{a+b}\lambda(s)e^{-\int_{0}^{s}\lambda(u)du}(s+M_{b}(s))ds
\\
&\geq
a\lambda(0)\int_{a}^{b}\lambda(s)e^{-\int_{0}^{s}\lambda(u)du}(s+M_{b}(s))ds
\\
&=a\lambda(0)\int_{0}^{b}\lambda(s)e^{-\int_{0}^{s}\lambda(u)du}(s+M_{b}(s))ds
-a\lambda(0)\int_{0}^{a}\lambda(s)e^{-\int_{0}^{s}\lambda(u)du}(s+M_{b}(s))ds
\\
&\geq
a\lambda(0)\int_{0}^{b}\lambda(s)e^{-\int_{0}^{s}\lambda(u)du}(s+M_{b}(s))ds
-a^{2}\lambda(0)\bar{\lambda}f_{+}\left(\frac{1}{2f_{+}\bar{\lambda}}+\frac{e^{\bar{\lambda}f_{+}b}-1}{\bar{\lambda}f_{+}}\right)
\\
&
=a\lambda(0)\int_{0}^{b}\lambda(s)e^{-\int_{0}^{s}\lambda(u)du}(s+M_{b}(s))ds
-a^{2}\lambda(0)\left(e^{\bar{\lambda}f_{+}b}-\frac{1}{2}\right).
\end{align*}
Next, we can compute that
\begin{align*}
&\int_{b}^{a+b}\lambda(s)e^{-\int_{0}^{s}\lambda(u)du}(s+M_{b}(s))ds
\\
&\geq
a\lambda(b)e^{-\int_{0}^{a+b}\lambda(u)du}(b+M_{b}(b))
\\
&=a\lambda(b)e^{-\int_{0}^{b}\lambda(u)du}(b+M_{b}(b))e^{-\int_{b}^{a+b}\lambda(u)du}
\\
&\geq
a\lambda(b)e^{-\int_{0}^{b}\lambda(u)du}(b+M_{b}(b))\left(1-\int_{b}^{a+b}\lambda(u)du\right)
\\
&\geq
a\lambda(b)e^{-\int_{0}^{b}\lambda(u)du}(b+M_{b}(b))
-a^{2}\lambda(b)e^{-\int_{0}^{b}\lambda(u)du}(b+M_{b}(b))\bar{\lambda}f_{+}
\\
&\geq
a\lambda(b)e^{-\int_{0}^{b}\lambda(u)du}(b+M_{b}(b))
-a^{2}\lambda(b)\left(b+\frac{e^{\bar{\lambda}f_{+}b}-1}{\bar{\lambda}f_{+}}\right)\bar{\lambda}f_{+}
\\
&=a\lambda(b)e^{-\int_{0}^{b}\lambda(u)du}(b+M_{b}(b))
-a^{2}\lambda(b)\left(b\bar{\lambda}f_{+}-1+e^{\bar{\lambda}f_{+}b}\right).
\end{align*}
Finally, we can compute that 
\begin{align*}
&-\int_{0}^{a}\lambda(s)e^{-\int_{0}^{s}\lambda(u)du}(s+M_{b}(s))ds
\\
&\geq
-a\lambda(0)(a+M_{b}(a))
=-a\lambda(0)M_{b}(0)+a\lambda(0)(M_{b}(0)-M_{b}(a))-a^{2}\lambda(0).
\end{align*}
Putting everything together, we get
\begin{align*}
M_{b}(a)-M_{b}(0)
&\geq
be^{-\int_{0}^{b}\lambda(s)ds}\left(a\lambda(b)-a\lambda(0)\right)
-be^{-\int_{0}^{b}\lambda(s)ds}a^{2}L
-a\left(1-e^{-\int_{0}^{b}\lambda(u)du}\right)
-a^{2}\bar{\lambda}f_{+}
\\
&\qquad
+a\lambda(0)\int_{0}^{b}\lambda(s)e^{-\int_{0}^{s}\lambda(u)du}(s+M_{b}(s))ds
-a^{2}\lambda(0)\left(e^{\bar{\lambda}f_{+}b}-\frac{1}{2}\right)
\\
&\qquad
+a\lambda(b)e^{-\int_{0}^{b}\lambda(u)du}(b+M_{b}(b))
-a^{2}\lambda(b)\left(b\bar{\lambda}f_{+}-1+e^{\bar{\lambda}f_{+}b}\right)
\\
&\qquad
-a\lambda(0)M_{b}(0)+a\lambda(0)(M_{b}(0)-M_{b}(a))-a^{2}\lambda(0).
\end{align*}
Since $a\leq\frac{1}{2\bar{\lambda}f_{+}}\leq\frac{1}{2\lambda(0)}$
and by using 
\begin{equation*}
M_{b}(0)=be^{-\int_{0}^{b}\lambda(s)ds}
+\int_{0}^{b}\lambda(s)e^{-\int_{0}^{s}\lambda(u)du}(s+M_{b}(s))ds,
\end{equation*}
we get that
\begin{align*}
\frac{3}{2}(M_{b}(a)-M_{b}(0))
&\geq
a\left(\lambda(b)e^{-\int_{0}^{b}\lambda(u)du}M_{b}(b)-1+e^{-\int_{0}^{b}\lambda(u)du}\right)
-a^{2}\bar{\lambda}f_{+}-a^{2}bL
\\
&\qquad
-a^{2}\lambda(0)\left(e^{\bar{\lambda}f_{+}b}+\frac{1}{2}\right)
-a^{2}\lambda(b)\left(b\bar{\lambda}f_{+}-1+e^{\bar{\lambda}f_{+}b}\right)
\\
&\geq 
ae^{-\int_{0}^{b}\lambda(u)du}\left(e^{\lambda(b)b}-e^{\int_{0}^{b}\lambda(u)du}\right)
-a^{2}\bar{\lambda}f_{+}-a^{2}bL
\\
&\qquad
-a^{2}\lambda(0)\left(e^{\bar{\lambda}f_{+}b}+\frac{1}{2}\right)
-a^{2}\lambda(b)\left(b\bar{\lambda}f_{+}-1+e^{\bar{\lambda}f_{+}b}\right),
\end{align*}
and therefore
\begin{align*}
M_{a,b}(0)-M_{b,a}(0)
&\geq
M_{b}(a)-M_{b}(0)-\bar{\lambda}f_{+}a^{2}
\\
&\geq
\frac{2}{3}ae^{-\int_{0}^{b}\lambda(u)du}\left(e^{\lambda(b)b}-e^{\int_{0}^{b}\lambda(u)du}\right)
-\frac{5}{3}a^{2}\bar{\lambda}f_{+}-\frac{2}{3}a^{2}bL
\\
&\qquad
-\frac{2}{3}a^{2}\lambda(0)\left(e^{\bar{\lambda}f_{+}b}+\frac{1}{2}\right)
-\frac{2}{3}a^{2}\lambda(b)\left(b\bar{\lambda}f_{+}-1+e^{\bar{\lambda}f_{+}b}\right)>0,
\end{align*}
provided that
\begin{equation}
a<(\mathcal{M}_{2})^{-1}e^{-\int_{0}^{b}\lambda(u)du}\left(e^{\lambda(b)b}-e^{\int_{0}^{b}\lambda(u)du}\right),
\end{equation}
where $\mathcal{M}_{2}$ is defined in \eqref{M2:defn}. 
This completes the proof of (ii).
\hfill $\Box$


\subsection{Proof of Theorem~\ref{thm:n}}

Before we proceed to the proof of Theorem~\ref{thm:n},
let us first provide some intuition.
By letting $t=0$ and $\lambda(t)=\bar{\lambda} f(t)$ in \eqref{eqn:a:1:n}, 
we get:
\begin{align}
&M_{1:n}(0)=A_{n}e^{-\int_{0}^{A_{n}}\bar{\lambda} f(s)ds}
\nonumber
\\
&\qquad\qquad
+\int_{0}^{A_{1}}\bar{\lambda} f(s)e^{-\int_{0}^{s}\bar{\lambda} f(u)du}(s+M_{1:n}(s))ds
+\int_{A_{1}}^{A_{2}}\bar{\lambda} f(s)e^{-\int_{0}^{s}\bar{\lambda} f(u)du}(s+M_{2:n}(s))ds
\nonumber
\\
&\qquad\qquad\qquad\qquad
+\cdots+\int_{A_{n-1}}^{A_{n}}\bar{\lambda} f(s)e^{-\int_{0}^{s}\bar{\lambda} f(u)du}(s+M_{n:n}(s))ds.
\label{F:a:1:a:n:eqn}
\end{align}
For any permutation $\pi$ of $\{1,2,\ldots,n\}$, 
we have $A_{n}=A_{\pi(n)}$.
Therefore, from \eqref{F:a:1:a:n:eqn}, we have
\begin{align}
M_{1:n}(0)-M_{\pi(1):\pi(n)}(0)
&=\int_{0}^{A_{1}}\bar{\lambda} f(s)e^{-\int_{0}^{s}\bar{\lambda} f(u)du}M_{1:n}(s)ds
+\int_{A_{1}}^{A_{2}}\bar{\lambda} f(s)e^{-\int_{0}^{s}\bar{\lambda} f(u)du}M_{2:n}(s)ds
\nonumber
\\
&\qquad
+\cdots+\int_{A_{n-1}}^{A_{n}}\bar{\lambda} f(s)e^{-\int_{0}^{s}\bar{\lambda} f(u)du}M_{n:n}(s)ds
\nonumber
\\
&\qquad
-\int_{0}^{A_{\pi(1)}}\bar{\lambda} f(s)e^{-\int_{0}^{s}\bar{\lambda} f(u)du}M_{\pi(1):\pi(n)}(s)ds
\nonumber
\\
&\qquad\qquad
-\int_{A_{\pi(1)}}^{A_{\pi(2)}}\bar{\lambda} f(s)e^{-\int_{0}^{s}\bar{\lambda} f(u)du}M_{\pi(2):\pi(n)}(s)ds
\nonumber
\\
&\qquad\qquad
-\cdots-\int_{A_{\pi(n-1)}}^{A_{\pi(n)}}\bar{\lambda} f(s)e^{-\int_{0}^{s}\bar{\lambda} f(u)du}M_{a_{\pi(n)}:a_{\pi(n)}}(s)ds.\label{F:diff:eqn}
\end{align}
Since there is no closed-form expression for $M_{i:n}(\cdot)$, 
it is not feasible to determine the sign of $M_{1:n}(0)-M_{\pi(1):\pi(n)}(0)$ directly. However, the key observation is that when $\bar{\lambda}$ is small, 
one can compute the first-order expansion of $M_{1:n}(0)-M_{\pi(1):\pi(n)}(0)$ in closed-form.
Indeed, as $\bar{\lambda}\rightarrow 0$, we have
$M_{i:n}(t)=A_{i:n}+o(1)$,
and $M_{\pi(i):\pi(n)}(t)=A_{\pi(i):\pi(n)}+o(1)$
uniformly for any $0\leq t\leq A_{n}$. 
Therefore, it follows from \eqref{F:diff:eqn} that we have
\begin{align}
&M_{1:n}(0)-M_{\pi(1):\pi(n)}(0)
\nonumber
\\
&=\bar{\lambda}A_{n}\int_{0}^{A_{1}}f(s)ds
+\bar{\lambda}A_{2:n}\int_{A_{1}}^{A_{2}}f(s)ds
+\cdots+\bar{\lambda} A_{n:n}\int_{A_{n-1}}^{A_{n}}f(s)ds
\nonumber
\\
&\qquad-\bar{\lambda}A_{\pi(n)}\int_{0}^{A_{\pi(1)}}f(s)ds
-\bar{\lambda}A_{\pi(2):\pi(n)}\int_{A_{\pi(1)}}^{A_{\pi(2)}}f(s)ds
-\cdots-\bar{\lambda}A_{\pi(n):\pi(n)}\int_{A_{\pi(n-1)}}^{A_{\pi(n)}}f(s)ds+o(\bar{\lambda})
\nonumber
\\
&=\bar{\lambda}A_{n}\int_{0}^{A_{1}}f(s)ds
+\bar{\lambda}A_{2:n}\left(\int_{0}^{A_{2}}f(s)ds-\int_{0}^{A_{1}}f(s)ds\right)
+\cdots+\bar{\lambda} A_{n:n}\left(\int_{0}^{A_{n}}f(s)ds-\int_{0}^{A_{n-1}}f(s)ds\right)
\nonumber
\\
&\qquad\qquad-\bar{\lambda}A_{\pi(n)}\int_{0}^{A_{\pi(1)}}f(s)ds
-\bar{\lambda}A_{\pi(2):\pi(n)}\left(\int_{0}^{A_{\pi(2)}}f(s)ds-\int_{0}^{A_{\pi(1)}}f(s)ds\right)
\nonumber
\\
&\qquad\qquad\qquad
-\cdots-\bar{\lambda} A_{\pi(n):\pi(n)}\left(\int_{0}^{A_{\pi(n)}}f(s)ds-\int_{0}^{A_{\pi(n-1)}}f(s)ds\right)
+o(\bar{\lambda})
\nonumber
\\
&=\bar{\lambda}\sum_{i=1}^{n}a_{i}\int_{0}^{A_{i}}f(s)ds
-\bar{\lambda}\sum_{i=1}^{n}a_{\pi(i)}\int_{0}^{A_{\pi(i)}}f(s)ds
+o(\bar{\lambda}).\label{F:diff:eqn:2}
\end{align}
We will make the above analysis rigorous in the proof of Theorem~\ref{thm:n} later on.
As a result, to determine the sign of $M_{1:n}(0)-M_{\pi(1):\pi(n)}(0)$
when $\bar{\lambda}$ is small, the key is to analyze the first term in the last line
in \eqref{F:diff:eqn:2}. Indeed, we have the following result.

\begin{lemma}\label{lem:ineq:n}
(i) If $f(t)$ is strictly decreasing for $t\leq A_{n}$, then
\begin{equation}
\sum_{i=1}^{n}a_{i}\int_{0}^{A_{i}}f(s)ds
<\sum_{i=1}^{n}a_{\pi(i)}\int_{0}^{A_{\pi(i)}}f(s)ds
\end{equation}
for any permutation $\pi$ such that $(\pi(1),\ldots,\pi(n))\neq(1,2,\ldots,n)$.

(ii) If $f(t)$ is strictly increasing for $t\leq A_{n}$, then
\begin{equation}
\sum_{i=1}^{n}a_{n+1-i}\int_{0}^{A_{n+1-i:n}}f(s)ds
<\sum_{i=1}^{n}a_{\pi(i)}\int_{0}^{A_{\pi(i)}}f(s)ds
\end{equation}
for any permutation $\pi$ such that $(\pi(1),\ldots,\pi(n))\neq(n,n-1,\ldots,1)$.
\end{lemma}

\noindent
{\bf Proof of Lemma~\ref{lem:ineq:n}}
Let us prove (i), 
and the argument to prove (ii) is similar.
To prove (i), it suffices to show that 
if there exists any $k$ such that $\pi(k)>\pi(k+1)$, 
then by defining $\tilde{\pi}$ so that 
$\tilde{\pi}(j)=\pi(j)$ for any $j\neq k,k+1$
and $\tilde{\pi}(k)=\pi(k+1)$
and $\tilde{\pi}(k+1)=\pi(k)$, we have
\begin{equation}\label{proof:key}
\sum_{i=1}^{n}a_{\tilde{\pi}(i)}\int_{0}^{A_{\tilde{\pi}(i)}}f(s)ds
<\sum_{i=1}^{n}a_{\pi(i)}\int_{0}^{A_{\pi(i)}}f(s)ds.
\end{equation}
Next, let us show that \eqref{proof:key} holds.
We can compute that
\begin{align*}
&\sum_{i=1}^{n}a_{\tilde{\pi}(i)}\int_{0}^{A_{\tilde{\pi}(i)}}f(s)ds
-\sum_{i=1}^{n}a_{\pi(i)}\int_{0}^{A_{\pi(i)}}f(s)ds
\\
&=a_{\tilde{\pi}(k)}\int_{0}^{A_{\tilde{\pi}(k)}}f(s)ds
+a_{\tilde{\pi}(k+1)}\int_{0}^{A_{\tilde{\pi}(k+1)}}f(s)ds
-a_{\pi(k)}\int_{0}^{A_{\pi(k)}}f(s)ds
-a_{\pi(k+1)}\int_{0}^{A_{\pi(k+1)}}f(s)ds
\\
&=a_{\pi(k+1)}\int_{0}^{A_{\pi(k-1)}+a_{\pi(k+1)}}f(s)ds
+a_{\pi(k)}\int_{0}^{A_{\pi(k+1)}}f(s)ds
\\
&\qquad\qquad\qquad\qquad
-a_{\pi(k)}\int_{0}^{A_{\pi(k)}}f(s)ds
-a_{\pi(k+1)}\int_{0}^{A_{\pi(k+1)}}f(s)ds
\\
&=a_{\pi(k)}\int_{A_{\pi(k)}}^{A_{\pi(k+1)}}f(s)ds
-a_{\pi(k+1)}\int_{A_{\pi(k-1)}+a_{\pi(k+1)}}^{A_{\pi(k+1)}}f(s)ds
\\
&=(a_{\pi(k)}-a_{\pi(k+1)})\int_{A_{\pi(k)}}^{A_{\pi(k+1)}}f(s)ds
-a_{\pi(k+1)}\int_{A_{\pi(k-1)}+a_{\pi(k+1)}}^{A_{\pi(k)}}f(s)ds.
\end{align*}
Since $f(t)$ is strictly decreasing in $t$ for $t\leq A_{n}=A_{\pi(n)}$,
and $\pi(k)>\pi(k+1)$, 
we have $A_{\pi(k)}>A_{\pi(k-1)}+a_{\pi(k+1)}$ and
\begin{align*}
&(a_{\pi(k)}-a_{\pi(k+1)})\int_{A_{\pi(k)}}^{A_{\pi(k+1)}}f(s)ds
-a_{\pi(k+1)}\int_{A_{\pi(k-1)}+a_{\pi(k+1)}}^{A_{\pi(k)}}f(s)ds
\\
&<(a_{\pi(k)}-a_{\pi(k+1)})\int_{A_{\pi(k)}}^{A_{\pi(k+1)}}f\left(A_{\pi(k)}\right)ds
-a_{\pi(k+1)}\int_{A_{\pi(k-1)}+a_{\pi(k+1)}}^{A_{\pi(k)}}f\left(A_{\pi(k)}\right)ds
\\
&=(a_{\pi(k)}-a_{\pi(k+1)})a_{\pi(k+1)}f\left(A_{\pi(k)}\right)
-a_{\pi(k+1)}(a_{\pi(k)}-a_{\pi(k+1)})f\left(A_{\pi(k)}\right)=0.
\end{align*}
Hence, we proved \eqref{proof:key}.
\hfill $\Box$

Before proceeding to the proof of Theorem~\ref{thm:n}, 
we state the following technical lemma.

\begin{lemma}\label{lem:bounds:n}
Assume $\bar{\lambda}\leq\frac{1}{2f_{+}a_{n}}$.
For any $t\geq 0$ and $i=1,2,\ldots,n$,
\begin{align}
A_{i:n}
\leq M_{i:n}(t)\leq A_{i:n}+\bar{\lambda} f_{+}(a_{i}^{2}+\cdots+a_{n}^{2}).
\end{align}
In particular, $M_{i:n}(t)\leq\frac{3}{2}A_{i:n}$.
\end{lemma}

\noindent
{\bf Proof of Lemma~\ref{lem:bounds:n}}
It follows from the definition of $M_{i:n}(t)$ any $t\geq 0$,
\begin{equation}
M_{i:n}(t)\geq a_{i}+a_{i+1}+\cdots+a_{n}=A_{i:n}.
\end{equation}
On the other hand, 
since $f_{-}\leq f(t)\leq f_{+}$
for any $t$, we have $\bar{\lambda} f_{-}\leq\lambda(t)\leq\bar{\lambda} f_{+}$
for any $t$, and as a result, for any $s\geq 0$, we have
\begin{align}
M_{i:n}(t)
&\leq  M_{a_{i};\bar{\lambda} f_{+}}+M_{a_{i+1};\bar{\lambda} f_{+}}
+\cdots+M_{a_{n};\bar{\lambda} f_{+}}
\nonumber
\\
&=\frac{1}{\bar{\lambda} f_{+}}(e^{\bar{\lambda} f_{+}a_{i}}-1)
+\frac{1}{\bar{\lambda} f_{+}}(e^{\bar{\lambda} f_{+}a_{i+1}}-1)
\cdots
+\frac{1}{\bar{\lambda} f_{+}}(e^{\bar{\lambda} f_{+}a_{n}}-1),\label{F:upper:bound}
\end{align}
where we applied Lemma~\ref{lem:const} to obtain the equality in equation~\eqref{F:upper:bound}.

Next, we notice that for any $0\leq x\leq\frac{1}{2}$, by Taylor's expansion,
\begin{equation}\label{e:x:ineq}
e^{x}-1=\sum_{k=1}^{\infty}\frac{x^{k}}{k!}
=x+\sum_{k=2}^{\infty}\frac{x^{k}}{k!}
\leq
x+\frac{1}{2}\sum_{k=2}^{\infty}x^{k}
=x+\frac{x^{2}}{2(1-x)}
\leq
x+x^{2}.
\end{equation}
Therefore, it follows from \eqref{F:upper:bound} and \eqref{e:x:ineq} that
\begin{align*}
M_{i:n}(t)
\leq
a_{i}+a_{i+1}+\cdots+a_{n}+\bar{\lambda} f_{+}(a_{i}^{2}+\cdots+a_{n}^{2})
\leq
\frac{3}{2}\left(a_{i}+a_{i+1}+\cdots+a_{n}\right)=\frac{3}{2}A_{i:n},
\end{align*}
provided that $\bar{\lambda}\leq\frac{1}{2f_{+}\max_{1\leq i\leq n}a_{i}}=\frac{1}{2f_{+}a_{n}}$.
This completes the proof.
\hfill $\Box$

Now we are ready to prove Theorem~\ref{thm:n}.

\noindent
{\bf Proof of Theorem~\ref{thm:n}}
(i) From equation~\eqref{F:diff:eqn}, using $e^{-x}\geq 1-x$ for any $x\geq 0$ and $f(t)\leq f_{+}$ for any $t$, we get
\begin{align}
M_{1:n}(0)-M_{\pi(1):\pi(n)}(0)
&\leq\int_{0}^{A_{1}}\bar{\lambda} f(s)M_{1:n}(s)ds
+\int_{A_{1}}^{A_{2}}\bar{\lambda} f(s)M_{2:n}(s)ds
+\cdots+\int_{A_{n-1}}^{A_{n}}\bar{\lambda} f(s)M_{n:n}(s)ds
\nonumber
\\
&\qquad
-\int_{0}^{A_{\pi(1)}}\bar{\lambda} f(s)\left(1-\int_{0}^{s}\bar{\lambda} f(u)du\right)M_{\pi(1):\pi(n)}(s)ds
\nonumber
\\
&\qquad\qquad
-\int_{A_{\pi(1)}}^{A_{\pi(2)}}\bar{\lambda} f(s)\left(1-\int_{0}^{s}\bar{\lambda} f(u)du\right)M_{\pi(2):\pi(n)}(s)ds
\nonumber
\\
&\qquad\qquad
-\cdots-\int_{A_{\pi(n-1)}}^{A_{\pi(n)}}\bar{\lambda} f(s)
\left(1-\int_{0}^{s}\bar{\lambda} f(u)du\right)M_{\pi(n):\pi(n)}(s)ds
\nonumber
\\
&=\int_{0}^{A_{1}}\bar{\lambda} f(s)M_{1:n}(s)ds
+\int_{A_{1}}^{A_{2}}\bar{\lambda} f(s)M_{2:n}(s)ds
+\cdots+\int_{A_{n-1}}^{A_{n}}\bar{\lambda} f(s)M_{n:n}(s)ds
\nonumber
\\
&\qquad
-\int_{0}^{A_{\pi(1)}}\bar{\lambda} f(s)M_{\pi(1):\pi(n)}(s)ds
-\int_{A_{\pi(1)}}^{A_{\pi(2)}}\bar{\lambda} f(s)M_{\pi(2):\pi(n)}(s)ds
\nonumber
\\
&\qquad\qquad
-\cdots-\int_{A_{\pi(n-1)}}^{A_{\pi(n)}}\bar{\lambda} f(s)M_{\pi(n):\pi(n)}(s)ds
\nonumber
\\
&\qquad
+\bar{\lambda}^{2}(f_{+})^{2}\int_{0}^{A_{\pi(1)}}sM_{\pi(1):\pi(n)}(s)ds
+\bar{\lambda}^{2}(f_{+})^{2}\int_{A_{\pi(1)}}^{A_{\pi(2)}}sM_{\pi(2):\pi(n)}(s)ds
\nonumber
\\
&\qquad\qquad
+\bar{\lambda}^{2}(f_{+})^{2}\int_{A_{\pi(n-1)}}^{A_{\pi(n)}}sM_{\pi(n):\pi(n)}(s)ds.
\nonumber
\end{align}
By Lemma~\ref{lem:bounds:n}, we get
\begin{align}
M_{1:n}(0)-M_{\pi(1):\pi(n)}(0)
&\leq\int_{0}^{A_{1}}\bar{\lambda} f(s)M_{1:n}(s)ds
+\int_{A_{1}}^{A_{2}}\bar{\lambda} f(s)M_{2:n}(s)ds
+\cdots+\int_{A_{n-1}}^{A_{n}}\bar{\lambda} f(s)M_{n:n}(s)ds
\nonumber
\\
&\qquad
-\int_{0}^{A_{\pi(1)}}\bar{\lambda} f(s)M_{\pi(1):\pi(n)}(s)ds
-\int_{A_{\pi(1)}}^{A_{\pi(2)}}\bar{\lambda} f(s)M_{\pi(2):\pi(n)}(s)ds
\nonumber
\\
&\qquad\qquad
-\cdots-\int_{A_{\pi(n-1)}}^{A_{\pi(n)}}\bar{\lambda} f(s)M_{\pi(n):\pi(n)}(s)ds
\nonumber
\\
&\qquad
+\bar{\lambda}^{2}(f_{+})^{2}\frac{3}{2}A_{\pi(n)}\int_{0}^{A_{\pi(1)}}sds
+\bar{\lambda}^{2}(f_{+})^{2}\frac{3}{2}A_{\pi(2):\pi(n)}\int_{A_{\pi(1)}}^{A_{\pi(2)}}sds
\nonumber
\\
&\qquad\qquad
+\bar{\lambda}^{2}(f_{+})^{2}\frac{3}{2}A_{\pi(n):\pi(n)}\int_{A_{\pi(n-1)}}^{A_{\pi(n)}}sds.\label{by:plug:1}
\end{align}
Note that $A_{\pi(n)}=A_{n}$,
and 
\begin{align}
&\bar{\lambda}^{2}(f_{+})^{2}\frac{3}{2}A_{\pi(n)}\int_{0}^{A_{\pi(1)}}sds
+\bar{\lambda}^{2}(f_{+})^{2}\frac{3}{2}A_{\pi(2):\pi(n)}\int_{A_{\pi(1)}}^{A_{\pi(2)}}sds
+\bar{\lambda}^{2}(f_{+})^{2}\frac{3}{2}A_{\pi(n):\pi(n)}\int_{A_{\pi(n-1)}}^{A_{\pi(n)}}sds\nonumber
\\
&\leq
\bar{\lambda}^{2}(f_{+})^{2}\frac{3}{2}
A_{\pi(n)}
\left(\int_{0}^{A_{\pi(1)}}sds
+\int_{A_{\pi(1)}}^{A_{\pi(2)}}sds
+\cdots+\int_{A_{\pi(n-1)}}^{A_{\pi(n)}}sds\right)
\nonumber
\\
&=\bar{\lambda}^{2}(f_{+})^{2}\frac{3}{2}
A_{n}\int_{0}^{A_{n}}sds
\nonumber
\\
&=\bar{\lambda}^{2}(f_{+})^{2}\frac{3}{4}(A_{n})^{3}.\label{by:plug:2}
\end{align}
Therefore, by plugging \eqref{by:plug:2} into \eqref{by:plug:1}, we have
\begin{align}
M_{1:n}(0)-M_{\pi(1):\pi(n)}(0)
&\leq\int_{0}^{A_{1}}\bar{\lambda} f(s)M_{1:n}(s)ds
+\int_{A_{1}}^{A_{2}}\bar{\lambda} f(s)M_{2:n}(s)ds
+\cdots+\int_{A_{n-1}}^{A_{n}}\bar{\lambda} f(s)M_{n:n}(s)ds
\nonumber
\\
&\qquad
-\int_{0}^{A_{\pi(1)}}\bar{\lambda} f(s)M_{\pi(1):\pi(n)}(s)ds
-\int_{A_{\pi(1)}}^{A_{\pi(2)}}\bar{\lambda} f(s)M_{\pi(2):\pi(n)}(s)ds
\nonumber
\\
&\qquad
-\cdots-\int_{A_{\pi(n-1)}}^{A_{\pi(n)}}\bar{\lambda} f(s)M_{\pi(n):\pi(n)}(s)ds
+\bar{\lambda}^{2}(f_{+})^{2}\frac{3}{4}(A_{n})^{3}.
\nonumber
\end{align}
By applying Lemma~\ref{lem:bounds:n}, we get
\begin{align}
&M_{1:n}(0)-M_{\pi(1):\pi(n)}(0)
\nonumber
\\
&\leq\int_{0}^{A_{1}}\bar{\lambda} f(s)\left(A_{n}+\bar{\lambda} f_{+}\sum_{i=1}^{n}a_{i}^{2}\right)ds
+\int_{A_{1}}^{A_{2}}\bar{\lambda} f(s)\left(A_{2:n}
+\bar{\lambda} f_{+}\sum_{i=2}^{n}a_{i}^{2}\right)ds
\nonumber
\\
&\qquad
+\cdots+\int_{A_{n-1}}^{A_{n}}\bar{\lambda} f(s)
\left(a_{n}+\bar{\lambda} f_{+}a_{n}^{2}\right)ds
\nonumber
\\
&\qquad
-\int_{0}^{A_{\pi(1)}}\bar{\lambda} f(s)A_{\pi(n)}ds
-\int_{A_{\pi(1)}}^{A_{\pi(2)}}\bar{\lambda} f(s)A_{\pi(2):\pi(n)}ds
-\cdots-\int_{A_{\pi(n-1)}}^{A_{\pi(n)}}\bar{\lambda} f(s)
a_{\pi(n)}ds
+\bar{\lambda}^{2}(f_{+})^{2}\frac{3}{4}(A_{n})^{3}
\nonumber
\\
&=
\bar{\lambda}\left(\sum_{i=1}^{n}a_{i}\int_{0}^{A_{i}}f(s)ds
-\sum_{i=1}^{n}a_{\pi(i)}\int_{0}^{A_{\pi(i)}}f(s)ds\right)
\nonumber
\\
&\qquad
+\int_{0}^{A_{1}}\bar{\lambda} f(s)\bar{\lambda} f_{+}\sum_{i=1}^{n}a_{i}^{2}ds
+\int_{A_{1}}^{A_{2}}\bar{\lambda} f(s)\bar{\lambda} f_{+}\sum_{i=2}^{n}a_{i}^{2}ds
+\cdots+\int_{A_{n-1}}^{A_{n}}\bar{\lambda} f(s)
\bar{\lambda} f_{+}a_{n}^{2}ds
+\bar{\lambda}^{2}(f_{+})^{2}\frac{3}{4}(A_{n})^{3}
\nonumber
\\
&\leq
\bar{\lambda}\left(\sum_{i=1}^{n}a_{i}\int_{0}^{A_{i}}f(s)ds
-\sum_{i=1}^{n}a_{\pi(i)}\int_{0}^{A_{\pi(i)}}f(s)ds\right)
\nonumber
\\
&\qquad
+\sum_{i=1}^{n}a_{i}^{2}\int_{0}^{A_{1}}\bar{\lambda} f(s)\bar{\lambda} f_{+}ds
+\sum_{i=1}^{n}a_{i}^{2}\int_{A_{1}}^{A_{2}}\bar{\lambda} f(s)\bar{\lambda} f_{+}ds
\nonumber
\\
&\qquad
+\cdots+\sum_{i=1}^{n}a_{i}^{2}\int_{A_{n-1}}^{A_{n}}\bar{\lambda} f(s)
\bar{\lambda} f_{+}ds
+\bar{\lambda}^{2}(f_{+})^{2}\frac{3}{4}(A_{n})^{3}
\nonumber
\\
&\leq
\bar{\lambda}\left(\sum_{i=1}^{n}a_{i}\int_{0}^{A_{i}}f(s)ds
-\sum_{i=1}^{n}a_{\pi(i)}\int_{0}^{A_{\pi(i)}}f(s)ds\right)
+\bar{\lambda}^{2}(f_{+})^{2}A_{n}\sum_{i=1}^{n}a_{i}^{2}
+\bar{\lambda}^{2}(f_{+})^{2}\frac{3}{4}\left(A_{n}\right)^{3}.
\nonumber
\end{align}
By Lemma~\ref{lem:ineq:n}, $\sum_{i=1}^{n}a_{i}\int_{0}^{A_{i}}f(s)ds
<\sum_{i=1}^{n}a_{\pi(i)}\int_{0}^{A_{\pi(i)}}f(s)ds$
and therefore $M_{1:n}(0)<M_{\pi(1):\pi(n)}(0)$
provided that
\begin{equation*}
\bar{\lambda}<\frac{
\sum_{i=1}^{n}a_{\pi(i)}\int_{0}^{A_{\pi(i)}}f(s)ds
-\sum_{i=1}^{n}a_{i}\int_{0}^{A_{i}}f(s)ds}
{(f_{+})^{2}\sum_{i=1}^{n}a_{i}^{2}\sum_{i=1}^{n}a_{i}
+(f_{+})^{2}\frac{3}{4}\left(\sum_{i=1}^{n}a_{i}\right)^{3}}.
\end{equation*}
This completes the proof of (i).

(ii) The proof is similar to part (i) and is hence omitted here.
\hfill $\Box$


\subsection{Proof of Theorem~\ref{thm:n:short}}

\noindent
{\bf Proof of Theorem~\ref{thm:n:short}}
Recall that $a_{i}=\epsilon\bar{a}_{i}$, where $\bar{a}_{1}<\bar{a}_{2}<\cdots<\bar{a}_{n}$.

(i) By Theorem~\ref{thm:n}, if $f(t)$ is strictly decreasing for $t\leq A_{n}$, then
$M_{1:n}(0)<M_{\pi(1):\pi(n)}(0)$,
for any permutation $\pi$ such that $(\pi(1),\ldots,\pi(n))\neq(1,2,\ldots,n)$
provided that $\bar{\lambda}\leq\frac{1}{2f_{+}\epsilon\bar{a}_{n}}$ and
\begin{equation*}
\bar{\lambda}<\frac{
\sum_{i=1}^{n}\bar{a}_{\pi(i)}\int_{0}^{\epsilon\bar{A}_{\pi(i)}}f(s)ds
-\sum_{i=1}^{n}\bar{a}_{i}\int_{0}^{\epsilon\bar{A}_{i}}f(s)ds}
{\epsilon^{2}((f_{+})^{2}\bar{A}_{n}\sum_{i=1}^{n}\bar{a}_{i}^{2}
+(f_{+})^{2}\frac{3}{4}\left(\bar{A}_{n}\right)^{3})}.
\end{equation*}
For any sufficiently small $\epsilon>0$, 
$\bar{\lambda}\leq\frac{1}{2f_{+}\epsilon\bar{a}_{n}}$.
Moreover, by L'H\^{o}pital's rule,
\begin{align*}
&\lim_{\epsilon\rightarrow 0}
\frac{
\sum_{i=1}^{n}\bar{a}_{\pi(i)}\int_{0}^{\epsilon\bar{A}_{\pi(i)}}f(s)ds
-\sum_{i=1}^{n}\bar{a}_{i}\int_{0}^{\epsilon\bar{A}_{i}}f(s)ds}
{\epsilon^{2}((f_{+})^{2}\bar{A}_{n}\sum_{i=1}^{n}\bar{a}_{i}^{2}
+(f_{+})^{2}\frac{3}{4}\left(\bar{A}_{n}\right)^{3})}
\\
&=\lim_{\epsilon\rightarrow 0}
\frac{\sum_{i=1}^{n}\bar{a}_{\pi(i)}\bar{A}_{\pi(i)}f(\epsilon\bar{A}_{\pi(i)})
-\sum_{i=1}^{n}\bar{a}_{i}\bar{A}_{i}f(\epsilon\bar{A}_{i})}
{2\epsilon((f_{+})^{2}\bar{A}_{n}\sum_{i=1}^{n}\bar{a}_{i}^{2}
+(f_{+})^{2}\frac{3}{4}\left(\bar{A}_{n}\right)^{3})}
\\
&=\frac{\sum_{i=1}^{n}\bar{a}_{\pi(i)}\left(\bar{A}_{\pi(i)}\right)^{2}f'(0)
-\sum_{i=1}^{n}\bar{a}_{i}\left(\bar{A}_{i}\right)^{2}f'(0)}
{2((f_{+})^{2}\bar{A}_{n}\sum_{i=1}^{n}\bar{a}_{i}^{2}
+(f_{+})^{2}\frac{3}{4}\left(\bar{A}_{n}\right)^{3})}.
\end{align*}
Hence, we conclude that if for any permutation $\pi$ with $(\pi(1),\ldots,\pi(n))\neq(1,2,\ldots,n)$ and
\begin{equation*}
|f'(0)|>\frac{2\bar{\lambda}\left((f_{+})^{2}\bar{A}_{n}\sum_{i=1}^{n}\bar{a}_{i}^{2}
+(f_{+})^{2}\frac{3}{4}\left(\bar{A}_{n}\right)^{3}\right)}{\sum_{i=1}^{n}\bar{a}_{i}\left(\bar{A}_{i}\right)^{2}-\sum_{i=1}^{n}\bar{a}_{\pi(i)}\left(\bar{A}_{\pi(i)}\right)^{2}},
\end{equation*}
then $M_{1:n}(0)<M_{\pi(1):\pi(n)}(0)$, which completes the proof of part (i).

(ii) The proof is similar to part (i) and is omitted here.
\hfill $\Box$


\subsection{Proof of Theorem~\ref{thm:equal:n}}

\noindent
{\bf Proof of Theorem~\ref{thm:equal:n}}
Recall that $a_{1}<a_{2}<\ldots<a_{n}$.
Then, we have
\begin{equation*}
M_{1:n}(0)=M_{a_{1}}(0)+\mathbb{E}[M_{a_{2}}(\tau_{a_{1}})]
+\mathbb{E}[M_{a_{3}}(\tau_{a_{2}})]
+\cdots+\mathbb{E}[M_{a_{n}}(\tau_{a_{n-1}})].
\end{equation*}
Since $\tau_{a_{i}}\geq a_{i}>t_{0}$ for any $i$, we have
\begin{equation*}
M_{1:n}(0)=M_{a_{1}}(0)+M_{a_{2};\lambda(t_{0})}
+M_{a_{3};\lambda(t_{0})}
+\cdots+M_{a_{n}:\lambda(t_{0})}.
\end{equation*}
Similarly, we can show that for any permutation $\pi$ of $\{1,\ldots,n\}$,
\begin{equation*}
M_{\pi(1):\pi(n)}(0)=M_{a_{\pi(1)}}(0)+M_{a_{\pi(2)};\lambda(t_{0})}
+M_{a_{\pi(3)};\lambda(t_{0})}
+\cdots+M_{a_{\pi(n)}:\lambda(t_{0})}.
\end{equation*}
If $\pi(1)=1$, then $M_{a_{1}}(0)=M_{a_{\pi(1)}}(0)$
and it follows that $M_{1:n}(0)=M_{\pi(1):\pi(n)}(0)$.
Otherwise, $a_{\pi(1)}>a_{1}$ since we assumed that $a_{1}<a_{2}<\ldots<a_{n}$.
Then, we can compute that
\begin{equation*}
M_{1:n}(0)-M_{\pi(1):\pi(n)}(0)
=M_{a_{1}}(0)-M_{a_{\pi(1)}}(0)+M_{a_{\pi(1)};\lambda(t_{0})}-M_{a_{1};\lambda(t_{0})}.
\end{equation*}
We claim that
\begin{equation}\label{claim:n:tasks}
M_{a_{1}}(0)+M_{a_{\pi(1)};\lambda(t_{0})}=M_{a_{\pi(1)}}(0)+M_{a_{1};\lambda(t_{0})}.
\end{equation}
This can be shown by the coupling argument. 
If $\pi(1)=1$, then $a_{\pi(1)}=a_{1}$ and \eqref{claim:n:tasks} trivially holds.
Otherwise, $a_{\pi(1)}>a_{1}$ since we assumed that $a_{1}<a_{2}<\ldots<a_{n}$.
Since $a_{\pi(1)}>a_{1}$, starting from time $0$, to the finish task of length $a_{\pi(1)}$,
it must cover the task of length $a_{1}$ first, and the expected value of time to finish 
the task of length $a_{\pi(1)}$ after finishing the task of length $a_{1}$ is given by $M_{a_{\pi(1)}}(0)-M_{a_{1}}(0)$. 
However, after the task of length $a_{1}$ is finished, the disruption rate is constant $\lambda(t_{0})$, 
thus $M_{a_{\pi(1)}}(0)-M_{a_{1}}(0)$ equals the expected value of time to finish 
the task of length $a_{\pi(1)}$ after finishing the task of length $a_{1}$ starting
from time $0$ with constant disruption rate $\lambda(t_{0})$.
Since $\mathbb{E}[M_{a_{\pi(1)}}(\tau_{a_{1}})]=M_{a_{\pi(1)};\lambda(t_{0})}$
and $\mathbb{E}[M_{a_{1}}(\tau_{a_{\pi(1)}})]=M_{a_{1};\lambda(t_{0})}$, 
$\mathbb{E}[M_{a_{\pi(1)}}(\tau_{a_{1}})]-\mathbb{E}[M_{a_{1}}(\tau_{a_{\pi(1)}})]$
equals the expected value of time to finish 
the task of length $a_{\pi(1)}$ after finishing the task of length $a_{1}$ starting
from time $0$ with constant disruption rate $\lambda(t_{0})$. Therefore equation~\eqref{claim:n:tasks} holds
and the proof is complete.
\hfill $\Box$


\subsection{Proof of Theorem~\ref{thm:one:break:n}}

\noindent
{\bf Proof of Theorem~\ref{thm:one:break:n}}
Let us recall that we have
\begin{align}
R_{1:n}(0)&=A_{n}e^{-\int_{0}^{A_{n}}\lambda(s)ds}
+\int_{0}^{A_{1}}\lambda(s)e^{-\int_{0}^{s}\lambda(u)du}(s+A_{n})ds
\nonumber
\\
&\quad
+\int_{A_{1}}^{A_{2}}\lambda(s)e^{-\int_{0}^{s}\lambda(u)du}(s+A_{2:n})ds
+\cdots+\int_{A_{n-1}}^{A_{n}}\lambda(s)e^{-\int_{0}^{s}\lambda(u)du}(s+A_{n:n})ds.
\nonumber
\end{align}
We can rewrite $R_{1:n}(0)$ as
\begin{align}
R_{1:n}(0)
&=A_{n}e^{-\int_{0}^{A_{n}}\lambda(s)ds}
+\int_{0}^{A_{1}}p(s)(s+A_{n})ds
+\int_{A_{1}}^{A_{2}}p(s)(s+A_{n})ds
-A_{1}\int_{A_{1}}^{A_{2}}p(s)ds
\nonumber
\\
&\qquad\qquad\qquad
+\cdots+\int_{A_{n-1}}^{A_{n}}p(s)(s+A_{n})ds
-\int_{A_{n-1}}^{A_{n}}p(s)A_{n-1}ds
\nonumber
\\
&=A_{n}e^{-\int_{0}^{A_{n}}\lambda(s)ds}
+\int_{0}^{A_{n}}p(s)(s+A_{n})ds
-A_{1}\int_{A_{1}}^{A_{2}}p(s)ds
-\cdots-A_{n-1}\int_{A_{n-1}}^{A_{n}}p(s)ds.
\nonumber
\end{align}
We can rewrite further $R_{1:n}(0)$ as
\begin{align}
R_{1:n}(0)
&=A_{n}e^{-\int_{0}^{A_{n}}\lambda(s)ds}
+\int_{0}^{A_{n}}p(s)(s+A_{n})ds
-A_{1}\left(\int_{0}^{A_{2}}p(s)ds-\int_{0}^{A_{1}}p(s)ds\right)
\nonumber
\\
&\qquad\qquad\qquad
-\cdots-A_{n-1}\left(\int_{0}^{A_{n}}p(s)ds-\int_{0}^{A_{n-1}}p(s)ds\right)
\nonumber
\\
&=A_{n}e^{-\int_{0}^{A_{n}}\lambda(s)ds}
+\int_{0}^{A_{n}}p(s)(s+A_{n})ds
+\sum_{i=1}^{n}a_{i}\int_{0}^{A_{i}}p(s)ds-A_{n}\int_{0}^{A_{n}}p(s)ds.
\nonumber
\end{align}
Therefore, together with the identity $A_{n}=A_{\pi(n)}$, 
we can compute that for any permutation $\pi$ of $\{1,2,\ldots,n\}$, we have
\begin{align*}
R_{1:n}(0)
-R_{\pi(1):\pi(n)}(0)
=\sum_{i=1}^{n}\left(a_{i}\int_{0}^{A_{i}}p(s)ds-a_{\pi(i)}\int_{0}^{A_{\pi(i)}}p(s)ds\right).
\end{align*}
It then follows from Lemma~\ref{lem:ineq:n}
that if $p(t)$ is strictly decreasing for $t\leq A_{n}$, then
$R_{1:n}(0)<R_{\pi(1):\pi(n)}(0)$,
for any permutation $\pi$ such that $(\pi(1),\ldots,\pi(n))\neq(1,2,\ldots,n)$,
and if $p(t)$ is strictly increasing for $t\leq A_{n}$, then 
$R_{n:1}(0)<R_{\pi(1):\pi(n)}(0)$,
for any permutation $\pi$ such that $(\pi(1),\ldots,\pi(n))\neq(n,n-1,\ldots,1)$.
The proof is complete.
\hfill $\Box$

\end{document}